\documentclass[10pt]{article}

\usepackage{geometry}
\usepackage{latexsym}
\usepackage{mathrsfs}
\usepackage{amsfonts}
\usepackage{amssymb}
\usepackage{subfigure}    
\usepackage{graphicx}
\usepackage{float}
\usepackage{multirow}
\usepackage{slashbox}
\usepackage{bm}
\usepackage{amsmath}
\usepackage{amsthm}
\usepackage{color}
\usepackage[subnum]{cases}
\usepackage{rotating}
\usepackage{enumitem}
\usepackage{accents}
\usepackage{algorithm}
\usepackage{algorithmic}

\usepackage[pdftex,
plainpages=false,
bookmarks,
bookmarksnumbered,
colorlinks=true,
linkcolor=black,
citecolor=black,
urlcolor=blue,
filecolor=black
]
{hyperref}

\makeatletter
\@addtoreset{equation}{section}
\makeatother

\def \beginproof{\par\noindent {\bf Proof.}\ }
\def \endproof{\hskip .5cm $\Box$ \vskip .5cm}

\makeatletter
\def \wideubar{\underaccent{{\cc@style\underline{\mskip15mu}}}}
\def \widebar{\accentset{{\cc@style\underline{\mskip11mu}}}}
\makeatother

\newcommand{\email}[1]{\protect\href{mailto:#1}{#1}}
\graphicspath{{images/}}

\begin{document}

\newtheorem{property}{Property}[section]
\newtheorem{proposition}{Proposition}[section]
\newtheorem{append}{Appendix}[section]
\newtheorem{definition}{Definition}[section]
\newtheorem{lemma}{Lemma}[section]
\newtheorem{corollary}{Corollary}[section]
\newtheorem{theorem}{Theorem}[section]
\newtheorem{remark}{Remark}[section]
\newtheorem{problem}{Problem}[section]
\newtheorem{assumption}{Assumption}[section]
\newtheorem{example}{Example}[section]

\title{A Non-monotone Alternating Updating Method for A Class of Matrix Factorization Problems}
\author{Lei Yang\footnotemark[1] \and Ting Kei Pong\footnotemark[2] \and Xiaojun Chen\footnotemark[2]}
\date{}
\maketitle

\renewcommand{\thefootnote}{\fnsymbol{footnote}}
\footnotetext[1]{Institute of Operations Research and Analytics, National University of Singapore, Singapore. (\email{orayl@nus.edu.sg}).}
\footnotetext[2]{Department of Applied Mathematics, The Hong Kong Polytechnic University, Hung Hom, Kowloon, Hong Kong, P.R. China. (\email{tk.pong@polyu.edu.hk}, \email{xiaojun.chen@polyu.edu.hk}). The second author's work was supported in part by Hong Kong Research Grants Council PolyU153085/16p. The third author's work was supported in part by Hong Kong Research Grants Council PolyU153000/15p.}
\renewcommand{\thefootnote}{\arabic{footnote}}

\begin{abstract}
In this paper we consider a general matrix factorization model which covers a large class of existing models with many applications in areas such as machine learning and imaging sciences. To solve this possibly nonconvex, nonsmooth and non-Lipschitz problem, we develop a non-monotone alternating updating method based on a potential function. Our method essentially updates two blocks of variables in turn by inexactly minimizing this potential function, and updates another auxiliary block of variables using an explicit formula. The special structure of our potential function allows us to take advantage of efficient computational strategies for non-negative matrix factorization to perform the alternating minimization over the two blocks of variables. A suitable line search criterion is also incorporated to improve the numerical performance. Under some mild conditions, we show that the line search criterion is well defined, and establish that the sequence generated is bounded and any cluster point of the sequence is a stationary point. Finally, we conduct some numerical experiments using real datasets to compare our method with some existing efficient methods for non-negative matrix factorization and matrix completion. The numerical results show that our method can outperform these methods for these specific applications.

\vspace{5mm}
\noindent {\bf Keywords:}~~Matrix factorization; non-monotone line search; stationary point;
alternating updating.

\end{abstract}


\section{Introduction}
\setcounter{section}{1}

In this paper we consider a class of matrix factorization problems, which can be modeled as
\begin{eqnarray}\label{MFPmodel}
\min\limits_{X,Y}~\mathcal{F}(X,Y):=\Psi(X)+\Phi(Y)+\frac{1}{2}\left\|\mathcal{A}(XY^{\top})-\bm{b}\right\|^2,
\end{eqnarray}
where $X\in\mathbb{R}^{m\times r}$ and $Y\in\mathbb{R}^{n\times r}$ are decision variables with $r \leq \min\{m,\,n\}$, the functions $\Psi: \mathbb{R}^{m\times r} \rightarrow \mathbb{R}\cup\{\infty\}$ and $\Phi: \mathbb{R}^{n\times r} \rightarrow \mathbb{R}\cup\{\infty\}$ are proper closed but possibly \textit{nonconvex}, \textit{nonsmooth} and \textit{non-Lipschitz}, $\bm{b}\in\mathbb{R}^{q}$ is a given vector and $\mathcal{A}: \mathbb{R}^{m \times n} \rightarrow \mathbb{R}^{q}$ is a linear map with $q \leq mn$ and $\mathcal{A}\mathcal{A}^*=\mathcal{I}_q$ ($\mathcal{I}_q$ denotes the identity map from $\mathbb{R}^{q}$ to $\mathbb{R}^{q}$). Model \eqref{MFPmodel} covers many existing widely-studied models in many application areas such as machine learning \cite{uhzb2016generalized} and imaging sciences \cite{zyzy2015low}. In particular, $\Psi(X)$ and $\Phi(Y)$ can be various regularizers for inducing desired structures, and $\mathcal{A}$ can be suitably chosen to model different scenarios. For example, when $\Psi(X)$ and $\Phi(Y)$ are chosen as the indicator functions (see the next section for notation and definitions) for $\mathcal{X}=\{X\in\mathbb{R}^{m\times r}:X \geq 0\}$ and $\mathcal{Y}=\{Y\in\mathbb{R}^{n\times r}:Y \geq 0\}$, respectively, and $\mathcal{A}$ is the identity map, \eqref{MFPmodel} reduces to the non-negative matrix factorization (NMF) problem, which has been widely used in data mining applications to provide interpretable decompositions of data. NMF was first introduced by Paatero and Tapper \cite{pt1994positive}, and then popularized by Lee and Seung \cite{ls1999learning}. The basic task of NMF is to find two nonnegative matrices $X\in\mathbb{R}^{m\times r}_{+}$ and $Y\in\mathbb{R}^{n\times r}_{+}$ such that $M \approx XY^{\top}$ for a given nonnegative data matrix $M\in\mathbb{R}^{m \times n}_{+}$. We refer readers to \cite{bblpp2007algorithms,g2011nonnegative,g2014the,ls2001algorithms,wz2013nonnegative} for more information on NMF and its variants. Another example of \eqref{MFPmodel} arises in recent models of the matrix completion (MC) problem (see \cite{slc2016scalable,slc2016tractable,slc2016unified}), where $\Psi(X)$ and $\Phi(Y)$ are chosen as the Schatten-$p_1$ quasi-norm and the Schatten-$p_2$ quasi-norm for suitable $p_1, p_2>0$, respectively, and $\mathcal{A}$ is the sampling map. The MC problem aims to recover an unknown low rank matrix from a sample of its entries and arises in various applications (see, for example, \cite{bltyw2006semidefinite,lv2009interior,rs2005fast,s2004learning}). Many widely-studied models for MC are based on nuclear-norm minimization \cite{cr2009exact,ct2010the,rfp2010guaranteed}, or, more generally, Schatten-$p$ ($0<p\leq1$) $\mathrm{(quasi-)norm}$ minimization \cite{lxy2013improved,mf2012iterative,zhz2013restricted}. Recently, models based on low-rank matrix factorization such as \eqref{MFPmodel} have become popular because singular value decompositions or eigenvalue decompositions of huge ($m\times n$) matrices are not required for solving these models (see, for example, \cite{kmo2010matrix,slc2016scalable,slc2016tractable,slc2016unified,sl2016guaranteed,wyz2012solving}). More examples of \eqref{MFPmodel} can be found in recent surveys \cite{uhzb2016generalized,zyzy2015low}.

Problem \eqref{MFPmodel} is in general nonconvex (even when $\Psi$, $\Phi$ are convex) and NP-hard\footnote{Problem \eqref{MFPmodel} is NP-hard because it contains NMF as a special case, which is NP-hard in general \cite{v2009on}.}. Therefore, in this paper, we focus on finding a \textit{stationary point} of the objective $\mathcal{F}$ in \eqref{MFPmodel}. Note that $\mathcal{F}$ involves two blocks of variables. This kind of structure has been widely studied in the literature; see, for example, \cite{abs2013convergence,bst2014proximal,h2016decomposing,jlmz2016structured,xy2013a,xywz2012an,z2010an}. One popular class of methods for tackling this kind of problems is the alternating direction method of multipliers (ADMM) (see, for example,  \cite{xywz2012an,z2010an}), in which each iteration consists of an alternating minimization of an augmented Lagrangian function that involves $X$, $Y$ and some auxiliary variables, followed by updates of the associated multipliers. However, the conditions presented in \cite{xywz2012an,z2010an} that guarantee convergence of the ADMM are too restrictive. Moreover, updating the auxiliary variables and the multipliers can be expensive for large-scale problems. Another class of methods for \eqref{MFPmodel} is the alternating-minimization-based (or block-coordinate-descent-type) methods (see \cite{abs2013convergence,bst2014proximal,cza2007hierarchical,gg2012accelerated,lz2009fastnmf,llwy2012sparse,xy2013a}), which alternately (exactly or inexactly) minimizes $\mathcal{F}(X,Y)$ over each block of variables and converges under some mild conditions. When $\mathcal{A}$ is not the identity map, the majorization technique can be used to simplify the subproblems. Some representative algorithms of this class are proximal alternating linearized minimization (PALM) \cite{bst2014proximal}, hierarchical alternating least squares (HALS) (for NMF only; see \cite{cza2007hierarchical,gg2012accelerated,lz2009fastnmf,llwy2012sparse}) and block coordinate descent (BCD) \cite{xy2013a}. Comparing with ADMM, it was reported in \cite{xy2013a} that BCD outperforms ADMM in both CPU time and solution quality for NMF.

PALM, HALS and BCD are currently the state-of-the-art algorithms for solving problems of the form \eqref{MFPmodel}. In this paper, we develop a new iterative method for \eqref{MFPmodel}, which, according to our numerical experiments in Section \ref{secnum}, outperforms HALS and BCD for NMF, and PALM for MC. Our method is based on the following potential function (specifically constructed for $\mathcal{F}$ in \eqref{MFPmodel}):
\begin{eqnarray}\label{defpofun}
\Theta_{\alpha,\beta} (X, Y, Z):=\Psi(X)+\Phi(Y)+\frac{\alpha}{2}\|XY^{\top}-Z\|_F^2
+\frac{\beta}{2}\left\|\mathcal{A}(Z)-\bm{b}\right\|^2,
\end{eqnarray}
where $\alpha$ and $\beta$ are real numbers. Instead of alternately (exactly or inexactly) minimizing $\mathcal{F}(X,Y)$ or the augmented Lagrangian function, our method alternately updates $X$ and $Y$ by inexactly minimizing $\Theta_{\alpha,\beta} (X, Y, Z)$ over $X$ and $Y$, and then updates $Z$ by an \textit{explicit formula}. Note that the coupled variables $XY^{\top}$ is now separated from $\mathcal{A}$ in our potential function. Thus, one can readily take advantage of efficient computational strategies for NMF, such as those used in HALS (see the ``hierarchical-prox" updating strategy in Section \ref{secalg}), for inexactly minimizing $\Theta_{\alpha,\beta} (X, Y, Z)$ over $X$ or $Y$. Furthermore, our method can be implemented for NMF and MC without explicitly forming the huge ($m\times n$) matrix $Z$ (see \eqref{substz1} and \eqref{substz2}) in each iteration. This significantly reduces the computational cost per iteration. Finally, a suitable non-monotone line search criterion, which is motivated by recent studies on non-monotone algorithms (see, for example, \cite{clp2016penalty,gzlhy2013a,wnf2009sparse}), is also incorporated to improve the numerical performance.

In the rest of this paper, we first present notation and preliminaries in Section \ref{secnot}. We then study the properties of our potential function $\Theta_{\alpha,\beta}$ in Section \ref{secpot}. Specifically, if $\mathcal{A}\mathcal{A}^*=\mathcal{I}_q$ and $\alpha$, $\beta$ are chosen such that $\alpha \mathcal{I} + \beta \mathcal{A}^*\mathcal{A} \succ 0$ and $\frac{1}{\alpha} + \frac{1}{\beta} = 1$, then the problem $\min\limits_{X,Y,Z} \left\{\Theta_{\alpha,\beta} (X, Y, Z)\right\}$ is equivalent to \eqref{MFPmodel} (see Theorem \ref{modeleq}). Furthermore, under the weaker conditions that $\mathcal{A}\mathcal{A}^*=\mathcal{I}_q$ and $\frac{1}{\alpha}+\frac{1}{\beta}=1$, we can show that (i) a stationary point of $\Theta_{\alpha,\beta}$ gives a stationary point of $\mathcal{F}$; (ii) a stationary point of $\mathcal{F}$ can be used to construct a stationary point of $\Theta_{\alpha,\beta}$ (see Theorem \ref{modelstaeq}). Thus, one can find a stationary point of $\mathcal{F}$ by finding a stationary point of $\Theta_{\alpha,\beta}$. In Section \ref{secalg}, we develop a non-monotone alternating updating method to find a stationary point of $\Theta_{\alpha,\beta}$, and hence of $\mathcal{F}$. The convergence analysis of our method is presented in Section \ref{secconv}. We show that our non-monotone line search criterion is well defined and any cluster point of the sequence generated by our method is a stationary point of $\mathcal{F}$ under some mild conditions. Section \ref{secnum} gives numerical experiments to evaluate the performance of our method for NMF and MC on real datasets. Our computational results illustrate the efficiency of our method. Finally, some concluding remarks are given in Section \ref{secconc}.

\section{Notation and preliminaries}\label{secnot}

In this paper, for a vector $\bm{x}\in\mathbb{R}^m$, $x_i$ denotes its $i$-th entry, $\|\bm{x}\|$ denotes the Euclidean norm of $\bm{x}$ and $\mathrm{Diag}(\bm{x})$ denotes the diagonal matrix whose $i$-th diagonal element is $x_i$. For a matrix $X \in \mathbb{R}^{m \times n}$, $x_{ij}$ denotes the $ij$-th entry of $X$, $\bm{x}_j$ denotes the $j$-th column of $X$ and $\mathrm{tr}(X)$ denotes the trace of $X$. The Schatten-$p$ (quasi-)norm ($0<p<\infty$) of $X$ is defined as $\|X\|_{S_p}=\left(\sum^{\min(m,n)}_{i=1}\varsigma^p_i(X)\right)^{\frac{1}{p}}$, where $\varsigma_i(X)$ is the $i$-th singular value of $X$. For $p=2$, the Schatten-$2$ norm reduces to the Frobenius norm $\|X\|_F$, and for $p=1$, the Schatten-$1$ norm reduces to the nuclear norm $\|X\|_*$. Moreover, the spectral norm is denoted by $\|X\|$, which is the largest singular value of $X$; and the $\ell_1$-norm and $\ell_p$-quasi-norm ($0 < p < 1$) of $X$ are given by $\|X\|_1:=\sum_{i=1}^{m}\sum_{j=1}^{n} |x_{ij}|$ and $\|X\|_p:=\left( \sum_{i=1}^{m}\sum_{j=1}^{n} |x_{ij}|^p \right)^{\frac{1}{p}}$, respectively. For two matrices $X$ and $Y$ of the same size, we denote their trace inner product by $\langle X,\,Y\rangle := \sum_{i=1}^m\sum_{j=1}^nx_{ij}y_{ij}$. We also use $X \leq Y$ (resp., $X \geq Y$) to denote $x_{ij}\leq y_{ij}$ (resp., $x_{ij}\geq y_{ij}$) for all $(i,j)$. Furthermore, for a linear map $\mathcal{A}:\mathbb{R}^{m \times n}\to \mathbb{R}^{q}$, $\mathcal{A}^*$ denotes the adjoint linear map and $\|\mathcal{A}\|$ denotes the induced operator norm of $\mathcal{A}$, i.e., $\|\mathcal{A}\|=\sup\{\|\mathcal{A}(X)\|:\|X\|_F \leq 1\}$. A linear self-map $\mathcal{T}$ is said to be symmetric if $\mathcal{T}=\mathcal{T}^*$. For a symmetric linear self-map $\mathcal{T}: \mathbb{R}^{m \times n} \to \mathbb{R}^{m\times n}$, we say that $\mathcal{T}$ is positive definite, denoted by $\mathcal{T} \succ 0$, if $\langle X, \,\mathcal{T}(X) \rangle > 0$ for all $X\neq0$. The identity map from $\mathbb{R}^{m \times n}$ to $\mathbb{R}^{m \times n}$ is denoted by $\mathcal{I}$ and the identity map from $\mathbb{R}^{q}$ to $\mathbb{R}^{q}$ is denoted by $\mathcal{I}_q$. Finally, for a nonempty closed set $\mathcal{C}\subseteq\mathbb{R}^{m \times n}$, its indicator function $\delta_{\mathcal{C}}$ is defined by
\begin{eqnarray*}
\delta_{\mathcal{C}}(X) = \left\{
\begin{array}{ll}
0 &\quad\mathrm{if}~X\in\mathcal{C}, \\
+\infty &\quad\mathrm{otherwise}.
\end{array}\right.
\end{eqnarray*}

For an extended-real-valued function $f: \mathbb{R}^{m \times n} \rightarrow [-\infty,\infty]$, we say that it is \textit{proper} if $f(X) > -\infty$ for all $X \in \mathbb{R}^{m \times n}$ and its domain ${\rm dom}f:=\{X\in \mathbb{R}^{m \times n} : f(X) < \infty\}$ is nonempty. A function $f:\mathbb{R}^{m \times n} \rightarrow [-\infty,\infty]$ is level-bounded \cite[Definition 1.8]{rw1998variational} if for every $\alpha\in\mathbb{R}$, the set $\{X \in \mathbb{R}^{m \times n}:f(X)\leq\alpha\}$ is bounded (possibly empty). For a proper function $f:\mathbb{R}^{m \times n} \rightarrow (-\infty,\infty]$, we use the notation $Y \xrightarrow{f} X$ to denote $Y \rightarrow X$ (i.e., $\|Y-X\|_F\to0$) and $f(Y) \rightarrow f(X)$. The (limiting) subdifferential \cite[Definition~8.3]{rw1998variational} of $f$ at $X\in \mathrm{dom}f$ used in this paper, denoted by $\partial f(X)$, is defined as
\begin{eqnarray*}
\partial f(X):=\left\{ D \in \mathbb{R}^{m\times n}: \exists\,X^k \xrightarrow{f} X~\mathrm{and}~D^k \rightarrow D ~\mathrm{with}~D^k \in \widehat{\partial} f(X^k) ~\mathrm{for~all}~k\right\},
\end{eqnarray*}
where $\widehat{\partial} f(\widetilde{Y})$ denotes the Fr\'{e}chet subdifferential of $f$ at $\widetilde{Y}\in \mathrm{dom}f$, which is the set of all $D \in \mathbb{R}^{m\times n}$ satisfying
\begin{eqnarray*}
\liminf\limits_{Y \neq \widetilde{Y}, Y \rightarrow \widetilde{Y}} \frac{f(Y)-f(\widetilde{Y})-\langle D, \,Y-\widetilde{Y}\rangle}{\|Y-\widetilde{Y}\|_F} \geq 0.
\end{eqnarray*}
From the above definition, we can easily observe (see, for example, \cite[Proposition~8.7]{rw1998variational}) that
\begin{eqnarray}\label{robust}
\left\{ D\in\mathbb{R}^{m\times n}: \exists\, X^k \xrightarrow{f} X,~ D^k \rightarrow D,~ D^k \in \partial f(X^k) \right\} \subseteq \partial f(X).
\end{eqnarray}
When $f$ is continuously differentiable or convex, the above subdifferential coincides with the classical concept of derivative or convex subdifferential of $f$; see, for example, \cite[Exercise~8.8]{rw1998variational} and \cite[Proposition~8.12]{rw1998variational}. In this paper, we say that $X^{*}$ is \textit{stationary point} of $f$ if $0 \in \partial f(X^*)$.

For a proper closed function $g:\mathbb{R}^{m}\rightarrow (-\infty,\infty]$, the proximal mapping $\mathrm{Prox}_{g}: \mathbb{R}^m \rightarrow \mathbb{R}^m$ of $g$ is defined by
$\mathrm{Prox}_{g}(\bm{z}):=\mathop{\mathrm{Argmin}}\limits_{\bm{x}\in\mathbb{R}^m}
\left\{g(\bm{x}) + \frac{1}{2}\|\bm{x} - \bm{z}\|^2\right\}$.
For any $\nu>0$, the matrix shrinkage operator $\mathcal{S}_{\nu}: \mathbb{R}^{m\times n} \rightarrow \mathbb{R}^{m\times n}$ is defined by
\begin{eqnarray*}
\mathcal{S}_{\nu}(X):=U\textrm{Diag}(\bar{\bm{s}})V^\top~\mbox{\rm with}~\bar{s}_i=
\left\{\begin{array}{ll}s_i-\nu,& \mbox{\rm if}\; s_i - \nu > 0, \\
0,& \mbox{\rm otherwise}, \end{array}\right.
\end{eqnarray*}
where $U\in\mathbb{R}^{m\times t}$, $\bm{s}\in\mathbb{R}^{t}_+$ and $V\in\mathbb{R}^{n \times t}$ are given by the singular value decomposition of $X$, i.e, $X = U\textrm{Diag}(\bm{s})V^\top$.

We now present two propositions, which will be useful for developing our method in Section \ref{secalg}.
\begin{proposition}\label{inverselemma}
Suppose that $\mathcal{A}\mathcal{A}^*=\mathcal{I}_q$ and $\alpha(\alpha+\beta)\neq0$. Then, $\alpha \mathcal{I} + \beta \mathcal{A}^*\mathcal{A}$ is invertible and its inverse is given by $\frac{1}{\alpha}\mathcal{I} - \frac{\beta}{\alpha(\alpha+\beta)}\mathcal{A}^*\mathcal{A}$.
\end{proposition}
\beginproof
It is easy to check that $\frac{1}{\alpha}\mathcal{I} - \frac{\beta}{\alpha(\alpha+\beta)}\mathcal{A}^*\mathcal{A}$ is well defined since $\alpha(\alpha+\beta)\neq 0$, and that
$\big{(}\alpha \mathcal{I} + \beta \mathcal{A}^*\mathcal{A}\big{)}\left(\frac{1}{\alpha}\mathcal{I} - \frac{\beta}{\alpha(\alpha+\beta)}\mathcal{A}^*\mathcal{A}\right)=\mathcal{I}$. This completes the proof.
\endproof

\begin{proposition}\label{propuse}
Let $\psi:\mathbb{R}^{m}\rightarrow (-\infty,\infty]$ and $\phi:\mathbb{R}^{n}\rightarrow (-\infty,\infty]$ be proper closed functions. Given $P,Q\in\mathbb{R}^{m\times n}$ and $\bm{a}\in\mathbb{R}^{n}$, $\bm{b}\in\mathbb{R}^{m}$ with $\|\bm{a}\|\neq 0$, $\|\bm{b}\|\neq 0$, the following statements hold.
\begin{itemize}
\item[{\rm (i)}] The problem $\min\limits_{\bm{x}\in\mathbb{R}^{m}}\left\{\psi(\bm{x}) + \frac{1}{2}\|\bm{x}\bm{a}^{\top}-P\|_F^2\right\}$
    is equivalent to
    $\min\limits_{\bm{x}\in\mathbb{R}^{m}}\left\{\psi(\bm{x}) + \frac{\|\bm{a}\|^2}{2}\left\|\bm{x}-P\bm{a}/\|\bm{a}\|^2\right\|^2\right\}$;

\item[{\rm (ii)}] The problem $\min\limits_{\bm{y}\in\mathbb{R}^{n}}\left\{\phi(\bm{y})
    +\frac{1}{2}\|\bm{b}\bm{y}^{\top}-Q\|_F^2\right\}$
    is equivalent to
    $\min\limits_{\bm{y}\in\mathbb{R}^{n}}\left\{\phi(\bm{y})+ \frac{\|\bm{b}\|^2}{2}\left\|\bm{y}-Q^{\top}\bm{b}/\|\bm{b}\|^2\right\|^2\right\}$.
\end{itemize}
\end{proposition}
\beginproof
Statement (i) can be easily proved by noticing that
\begin{eqnarray*}
\begin{aligned}
\|\bm{x}\bm{a}^{\top}-P\|_F^2&=\|\bm{x}\bm{a}^{\top}\|_F^2-2\langle\bm{x}\bm{a}^{\top}, P\rangle+\|P\|_F^2
=\|\bm{a}\|^2\|\bm{x}\|^2-2\langle\bm{x}, P\bm{a}\rangle+\|P\|_F^2 \\
&=\|\bm{a}\|^2\left\|\bm{x}-P\bm{a}/\|\bm{a}\|^2\right\|^2
-\|P\bm{a}\|^2/\|\bm{a}\|^2+\|P\|^2_F.
\end{aligned}
\end{eqnarray*}
Then, statement (ii) can be easily proved by using statement (i) and $\|\bm{b}\bm{y}^{\top}-Q\|_F^2=\|\bm{y}\bm{b}^{\top}-Q^{\top}\|_F^2$.
\endproof

Before ending this section, we discuss the first-order necessary conditions for \eqref{MFPmodel}. First, from \cite[Exercise~8.8]{rw1998variational} and \cite[Proposition~10.5]{rw1998variational}, we see that
\begin{eqnarray*}
\partial \mathcal{F}(X,\,Y) = \left(
\begin{array}{c}
\partial \Psi(X)+\mathcal{A}^*\left(\mathcal{A}(XY^{\top})-\bm{b}\right)Y \\
\partial \Phi(Y)+\left(\mathcal{A}^*\left(\mathcal{A}(XY^{\top})-\bm{b}\right)\right)^{\top}X
\end{array}\right).
\end{eqnarray*}
Then, it follows from the generalized Fermat's rule \cite[Theorem~10.1]{rw1998variational} that any local minimizer $(\widebar{X}, \widebar{Y})$ of \eqref{MFPmodel} satisfies $0 \in \partial \mathcal{F}(\widebar{X}, \widebar{Y})$, i.e.,
\begin{eqnarray}\label{optcond}
\left\{\begin{aligned}
&0\in \partial \Psi(\widebar{X})+\mathcal{A}^*(\mathcal{A}(\widebar{X}\widebar{Y}^{\top})-\bm{b})\widebar{Y},\\
&0\in \partial \Phi(\widebar{Y})+(\mathcal{A}^*(\mathcal{A}(\widebar{X}\widebar{Y}^{\top})
-\bm{b}))^{\top}\widebar{X},
\end{aligned}\right.
\end{eqnarray}
which implies that $(\widebar{X}, \widebar{Y})$ is a stationary point of $\mathcal{F}$. In this paper, we focus on finding a stationary point $(X^*, Y^*)$ of $\mathcal{F}$, i.e., $(X^*, Y^*)$ satisfies \eqref{optcond} in place of $(\widebar{X}, \widebar{Y})$.

\section{The potential function for $\mathcal{F}$}\label{secpot}

In this section, we analyze the relation between $\mathcal{F}$ and its potential function $\Theta_{\alpha,\beta}$ defined in \eqref{defpofun}. Intuitively, $\Theta_{\alpha,\beta}$ originates from $\mathcal{F}$ by separating the coupled variables $XY^{\top}$ from the linear mapping $\mathcal{A}$ via introducing an auxiliary variable $Z$ and penalizing $XY^{\top}=Z$. We will see later that the stationary point of $\mathcal{F}$ can be characterized by the stationary point of $\Theta_{\alpha,\beta}$. Before proceeding, we prove the following technical lemma.

\begin{lemma}\label{supplem}
Suppose that $\mathcal{A}\mathcal{A}^*=\mathcal{I}_q$ and $\frac{1}{\alpha}+\frac{1}{\beta}=1$. Then, for any $(X, Y, Z)$ satisfying
\begin{eqnarray}\label{tildeZ}
Z = {\textstyle\left(\mathcal{I}-\frac{\beta}{\alpha+\beta}\mathcal{A}^*\mathcal{A}\right)}\left(XY^{\top}\right) + {\textstyle\frac{\beta}{\alpha+\beta}}\mathcal{A}^*(\bm{b}),
\end{eqnarray}
we have $\mathcal{F}(X,Y)=\Theta_{\alpha,\beta}(X,Y,Z)$.
\end{lemma}
\beginproof
First, from \eqref{tildeZ}, we have
\begin{eqnarray}
XY^{\top} - Z &=& {\textstyle\frac{\beta}{\alpha+\beta}}\mathcal{A}^*(\mathcal{A}(XY^{\top})-\bm{b}),  \label{extraeq1}\\
\mathcal{A}(Z)-\bm{b} &=& \mathcal{A}\left(XY^{\top} - {\textstyle\frac{\beta}{\alpha+\beta}}\mathcal{A}^*\mathcal{A}(XY^{\top})
+{\textstyle\frac{\beta}{\alpha+\beta}}\mathcal{A}^*(\bm{b})\right)-\bm{b}    \nonumber\\
&=&\mathcal{A}(XY^{\top}) - {\textstyle\frac{\beta}{\alpha+\beta}}\mathcal{A}\mathcal{A}^*\mathcal{A}(XY^{\top})
+{\textstyle\frac{\beta}{\alpha+\beta}}\mathcal{A}\mathcal{A}^*(\bm{b}) - \bm{b}   \nonumber\\
&=&{\textstyle\frac{\alpha}{\alpha+\beta}}\left(\mathcal{A}(XY^{\top}) - \bm{b}\right),    \label{extraeq2}
\end{eqnarray}
where the last equality follows from $\mathcal{A}\mathcal{A}^*=\mathcal{I}_q$. Then, we see that
\begin{eqnarray*}
\begin{aligned}
&~~\frac{\alpha}{2}\|XY^{\top}-Z\|_F^2
+\frac{\beta}{2}\left\|\mathcal{A}(Z)-\bm{b}\right\|^2 \\
&=\frac{\alpha}{2}\left\|{\textstyle\frac{\beta}{\alpha+\beta}}
\mathcal{A}^*(\mathcal{A}(XY^{\top})-\bm{b})\right\|_F^2
+\frac{\beta}{2}\left\|{\textstyle\frac{\alpha}{\alpha+\beta}}\left(\mathcal{A}(XY^{\top}) - \bm{b}\right)\right\|^2  \\
&={\textstyle\frac{\alpha\beta^2}{(\alpha+\beta)^2}}\cdot\frac{1}{2}
\left\|\mathcal{A}^*(\mathcal{A}(XY^{\top})-\bm{b})\right\|_F^2
+{\textstyle\frac{\alpha^2\beta}{(\alpha+\beta)^2}}\cdot\frac{1}{2}\left\| \mathcal{A}(XY^{\top}) - \bm{b}\right\|^2 \\
&={\textstyle\frac{\alpha\beta^2}{(\alpha+\beta)^2}}\cdot
\frac{1}{2}\left\|\mathcal{A}(XY^{\top})-\bm{b}\right\|^2
+{\textstyle\frac{\alpha^2\beta}{(\alpha+\beta)^2}}\cdot\frac{1}{2}\left\| \mathcal{A}(XY^{\top}) - \bm{b}\right\|^2 \\
&={\textstyle\frac{\alpha\beta}{\alpha+\beta}}\cdot\frac{1}{2}\left\|\mathcal{A}(XY^{\top})
-\bm{b}\right\|^2,
\end{aligned}
\end{eqnarray*}
where the first equality follows from \eqref{extraeq1} and \eqref{extraeq2}; and the third equality follows from $\mathcal{A}\mathcal{A}^*=\mathcal{I}_q$. This, together with $\frac{1}{\alpha}+\frac{1}{\beta}=1$ and the definitions of $\mathcal{F}$ and $\Theta_{\alpha,\beta}$ completes the proof.
\endproof

Based on the above lemma, we now establish the following property of $\Theta_{\alpha,\beta}$.

\begin{theorem}\label{modeleq}
Suppose that $\mathcal{A}\mathcal{A}^*=\mathcal{I}_q$. If $\alpha$ and $\beta$ are chosen such that $\alpha \mathcal{I} + \beta \mathcal{A}^*\mathcal{A} \succ 0$ and $\frac{1}{\alpha} + \frac{1}{\beta} = 1$, then the problem $\min\limits_{X,Y,Z}
\left\{\Theta_{\alpha,\beta} (X, Y, Z)\right\}$ is equivalent to \eqref{MFPmodel}.
\end{theorem}
\beginproof
First, it is easy to see from $\alpha \mathcal{I} + \beta \mathcal{A}^*\mathcal{A} \succ 0$ that the function $Z \longmapsto \Theta_{\alpha,\beta} (X, Y, Z)$ is strongly convex. Thus, for any fixed $X$ and $Y$, the optimal solution $Z^*$ to the problem $\min\limits_{Z}
\left\{\Theta_{\alpha,\beta} (X, Y, Z)\right\}$ exists and is unique, and can be obtained explicitly. Indeed, from the optimality condition, we have
\begin{eqnarray*}
\alpha(Z^*-XY^{\top}) + \beta\mathcal{A}^*(\mathcal{A}(Z^*)-\bm{b}) = 0.
\end{eqnarray*}
Then, since $\alpha \mathcal{I} + \beta \mathcal{A}^*\mathcal{A}$ is invertible (as $\alpha \mathcal{I}+\beta \mathcal{A}^*\mathcal{A}\succ0$), we see that
\begin{eqnarray*}
\begin{aligned}
Z^* &= \left(\alpha\mathcal{I}+\beta\mathcal{A}^*\mathcal{A}\right)^{-1}\left[\alpha XY^{\top}+\beta\mathcal{A}^*(\bm{b})\right] \\
&={\textstyle\left[\frac{1}{\alpha}\mathcal{I}-\frac{\beta}{\alpha(\alpha+\beta)}\mathcal{A}^*\mathcal{A}\right]}
\left[\alpha XY^{\top}+\beta\mathcal{A}^*(\bm{b})\right] \\
&={\textstyle\left(\mathcal{I}-\frac{\beta}{\alpha+\beta}\mathcal{A}^*\mathcal{A}\right)}(XY^{\top})
+\left[{\textstyle\frac{\beta}{\alpha}}\mathcal{A}^*(\bm{b})
-{\textstyle\frac{\beta^2}{\alpha(\alpha+\beta)}}\mathcal{A}^*\mathcal{A}\mathcal{A}^*(\bm{b})\right] \\
&={\textstyle\left(\mathcal{I}-\frac{\beta}{\alpha+\beta}\mathcal{A}^*\mathcal{A}\right)}(XY^{\top})
+{\textstyle\left[\frac{\beta}{\alpha}-\frac{\beta^2}{\alpha(\alpha+\beta)}\right]}\mathcal{A}^*(\bm{b})\\
&={\textstyle\left(\mathcal{I}-\frac{\beta}{\alpha+\beta}\mathcal{A}^*\mathcal{A}\right)}(XY^{\top})
+{\textstyle\frac{\beta}{\alpha+\beta}}\mathcal{A}^*(\bm{b}),
\end{aligned}
\end{eqnarray*}
where the second equality follows from Proposition \ref{inverselemma} and the fourth equality follows from $\mathcal{A}\mathcal{A}^*=\mathcal{I}_q$. This, together with Lemma \ref{supplem}, implies that $\mathcal{F}(X,Y)=\Theta_{\alpha,\beta}(X,Y,Z^*)$. Then, we have that
\begin{eqnarray*}
\min\limits_{X,Y,Z}
\left\{\Theta_{\alpha,\beta} (X, Y, Z)\right\}  =\min\limits_{X,Y}\left\{\min\limits_{Z}\left\{\Theta_{\alpha,\beta} (X, Y, Z)\right\}\right\}
=\min\limits_{X,Y}\left\{\Theta_{\alpha,\beta} (X, Y, Z^*)\right\}
=\min\limits_{X,Y}\left\{\mathcal{F}(X,Y)\right\}.
\end{eqnarray*}
This completes the proof.
\endproof

\begin{remark}\label{remark1}
From the proof of Lemma \ref{supplem}, we see that if $\Phi$ and $\Psi$ are the indicator functions of some nonempty closed sets, then $\mathcal{F}(X,Y)=\left(\frac{1}{\alpha}+\frac{1}{\beta}\right)\Theta_{\alpha,\beta}(X,Y,Z)$ holds with the special choice of $Z$ in \eqref{tildeZ} whenever $\mathcal{A}\mathcal{A}^*=\mathcal{I}_q$ and $\frac{1}{\alpha}+\frac{1}{\beta}>0$. Thus, the result in Theorem \ref{modeleq} remains valid whenever $\mathcal{A}\mathcal{A}^*=\mathcal{I}_q$ and $\alpha$, $\beta$ are chosen such that $\alpha \mathcal{I} + \beta \mathcal{A}^*\mathcal{A} \succ 0$ and $\frac{1}{\alpha} + \frac{1}{\beta} > 0$.
\end{remark}

It can be seen from Theorem \ref{modeleq} that \eqref{MFPmodel} is equivalent to minimizing $\Theta_{\alpha,\beta}$ with some suitable choices of $\alpha$ and $\beta$. On the other hand, we can also characterize the relation between the stationary points of $\mathcal{F}$ and $\Theta_{\alpha,\beta}$ under weaker conditions on $\alpha$ and $\beta$.

\begin{theorem}\label{modelstaeq}
Suppose that $\mathcal{A}\mathcal{A}^*=\mathcal{I}_q$ and $\alpha$, $\beta$ are chosen such that $\frac{1}{\alpha} + \frac{1}{\beta} = 1$. Then, the following statements hold.
\begin{itemize}
\item[{\rm (i)}] If $(X^*, Y^*, Z^*)$ is a stationary point of $\Theta_{\alpha,\beta}$, then $(X^*, Y^*)$ is a stationary point of $\mathcal{F}$;

\item[{\rm (ii)}] If $(X^*, Y^*)$ is a stationary point of $\mathcal{F}$, then $(X^*, Y^*, Z^*)$ is a stationary point of $\Theta_{\alpha,\beta}$, where $Z^*$ is given by
    \begin{eqnarray}\label{addZ}
    Z^* = {\textstyle\left(\mathcal{I}-\frac{\beta}{\alpha+\beta}\mathcal{A}^*\mathcal{A}\right)}\left(X^*(Y^*)^{\top}\right) + {\textstyle\frac{\beta}{\alpha+\beta}}\mathcal{A}^*(\bm{b}).
    \end{eqnarray}
\end{itemize}
\end{theorem}
\beginproof
First, if $(X^*, Y^*, Z^*)$ is a stationary point of $\Theta_{\alpha,\beta}$, then we have $0 \in \partial \Theta_{\alpha,\beta}(X^*, Y^*, Z^*)$, i.e.,
\begin{numcases}{}
0 \in \partial\Psi(X^*) + \alpha(X^*(Y^*)^{\top}-Z^*)Y^*,   \label{optcx}\\
0 \in \partial\Phi(Y^*) + \alpha(X^*(Y^*)^{\top}-Z^*)^{\top}X^*,  \label{optcy}\\
0 = \alpha(Z^*-X^*(Y^*)^{\top}) + \beta\mathcal{A}^*(\mathcal{A}(Z^*)-\bm{b}).  \label{optcz}
\end{numcases}
Since $\frac{1}{\alpha}+\frac{1}{\beta}=1$, we have $\alpha(\alpha+\beta)\neq0$ and hence $\alpha \mathcal{I} + \beta \mathcal{A}^*\mathcal{A}$ is invertible from Lemma \ref{inverselemma}. Then, using the same arguments in the proof of Theorem \ref{modeleq}, we see from \eqref{optcz} that $(X^*, Y^*, Z^*)$ satisfies \eqref{addZ}. Moreover, using \eqref{addZ} and the same arguments in \eqref{extraeq1} and \eqref{extraeq2}, we have
\begin{eqnarray}
X^*(Y^*)^{\top} - Z^* &=& {\textstyle\frac{\beta}{\alpha+\beta}}\mathcal{A}^*(\mathcal{A}(X^*(Y^*)^{\top})-\bm{b}),  \label{addeq1}\\
\mathcal{A}(Z^*)-\bm{b}  &=&{\textstyle\frac{\alpha}{\alpha+\beta}}\left(\mathcal{A}(X^*(Y^*)^{\top}) - \bm{b}\right).    \label{addeq2}
\end{eqnarray}
Thus, substituting \eqref{addeq1} into \eqref{optcx} and \eqref{optcy}, we see that
\begin{eqnarray}\label{reoptc}
\left\{\begin{aligned}
&0\in \partial \Psi(X^*)+{\textstyle\frac{\alpha\beta}{\alpha+\beta}}\mathcal{A}^*(\mathcal{A}(X^*(Y^*)^{\top}) - \bm{b})Y^*,\\
&0\in \partial \Phi(Y^*)+{\textstyle\frac{\alpha\beta}{\alpha+\beta}}\left(\mathcal{A}^*(\mathcal{A}(X^*(Y^*)^{\top}) - \bm{b})\right)^{\top}X^*.
\end{aligned}\right.
\end{eqnarray}
This together with $\frac{1}{\alpha} + \frac{1}{\beta} = 1$ implies $(X^*, Y^*)$ is a stationary point of $\mathcal{F}$. This proves statement (i).

We now prove statement (ii). First, if $(X^*, Y^*)$ is a stationary point of $\mathcal{F}$, then invoking $\frac{1}{\alpha}+\frac{1}{\beta}=1$ and \eqref{optcond}, we have \eqref{reoptc}. Next, we consider $(X^*, Y^*, Z^*)$ with $Z^*$ given by \eqref{addZ}. Then, $(X^*, Y^*, Z^*)$ satisfies \eqref{addeq1} and \eqref{addeq2}. Thus, substituting \eqref{addeq1} into \eqref{reoptc}, we obtain \eqref{optcx} and \eqref{optcy}. Moreover, we have from \eqref{addeq1} and \eqref{addeq2} that
\begin{eqnarray}\label{zequa}
\begin{aligned}
&\quad\alpha(Z^*-X^*(Y^*)^{\top}) + \beta\mathcal{A}^*(\mathcal{A}(Z^*)-\bm{b})  \\
&=-{\textstyle\frac{\alpha\beta}{\alpha+\beta}}\mathcal{A}^*\left((\mathcal{A}(X^*(Y^*)^{\top})-\bm{b}\right) + \beta \mathcal{A}^*\left({\textstyle\frac{\alpha}{\alpha+\beta}}\left(\mathcal{A}(X^*(Y^*)^{\top}) - \bm{b}\right)\right) = 0.
\end{aligned}
\end{eqnarray}
This together with \eqref{optcx} and \eqref{optcy} implies that $(X^*, Y^*, Z^*)$ is a stationary point of $\Theta_{\alpha,\beta}$. This proves statement (ii).
\endproof

\begin{remark}\label{conecond}
From the proof of Theorem \ref{modelstaeq}, one can see that if $\partial\Psi$ and $\partial\Phi$ are cones, Theorem \ref{modelstaeq} remains valid under the weaker conditions that $\mathcal{A}\mathcal{A}^*=\mathcal{I}_q$ and $\frac{1}{\alpha} + \frac{1}{\beta} > 0$.
\end{remark}

From Theorem \ref{modelstaeq}, we see that a stationary point of $\mathcal{F}$ can be obtained from a stationary point of $\Theta_{\alpha,\beta}$ with a suitable choice of $\alpha$ and $\beta$, i.e., $\frac{1}{\alpha} + \frac{1}{\beta} = 1$. Since the linear map $\mathcal{A}$ is no longer associated with the coupled variables $XY^{\top}$ in $\Theta_{\alpha,\beta}$, finding a stationary point of $\Theta_{\alpha,\beta}$ is conceivably easier. Thus, one can consider finding a stationary point of $\Theta_{\alpha,\beta}$ in order to find a stationary point of $\mathcal{F}$. Note that some existing alternating-minimization-based methods (see, for example, \cite{abs2013convergence,xy2013a}) can be used to find a stationary point of $\Theta_{\alpha,\beta}$, and hence of $\mathcal{F}$, under the conditions that $\mathcal{A}\mathcal{A}^*=\mathcal{I}_q$ and $\alpha$, $\beta$ are chosen so that $\alpha \mathcal{I} + \beta \mathcal{A}^*\mathcal{A} \succ 0$ and $\frac{1}{\alpha} + \frac{1}{\beta} = 1$. These conditions further imply that $\alpha > 1$ and $\beta = \frac{\alpha}{\alpha-1}> 1$. However, as we will see from our numerical results in Section \ref{secnum}, finding a stationary point of $\Theta_{\alpha,\beta}$ with $\alpha > 1$ can be slow. In view of this, in the next section, we develop a new non-monotone alternating updating method for finding a stationary of $\Theta_{\alpha,\beta}$ (and hence of $\mathcal{F}$) under the weaker conditions that $\mathcal{A}\mathcal{A}^*=\mathcal{I}_q$ and $\frac{1}{\alpha} + \frac{1}{\beta} = 1$. This allows more flexibilities in choosing $\alpha$ and $\beta$.

\section{Non-monotone alternating updating method}\label{secalg}

In this section, we consider a non-monotone alternating updating method (NAUM) for finding a stationary point of $\Theta_{\alpha,\beta}$ with $\frac{1}{\alpha} + \frac{1}{\beta} = 1$. Compared to existing alternating-minimization-based methods \cite{abs2013convergence,xy2013a} applied to $\Theta_{\alpha,\beta}$, which update $X$, $Y$, $Z$ by alternately solving subproblems related to $\Theta_{\alpha,\beta}$, NAUM updates $Z$ by an \textit{explicit formula} (see \eqref{Zupdate}) and updates $X$, $Y$ by solving subproblems related to $\Theta_{\alpha,\beta}$ in a Gauss-Seidel manner. Before presenting the complete algorithm, we first comment on the updates of $X$ and $Y$.

Let $(X^k, Y^k)$ denote the value of $(X, Y)$ after the $(k\mathrm{-}1)^{\mathrm{th}}$ iteration, and let $(U, V)$ denote the candidate for $(X^{k+1}, Y^{k+1})$ at the $k$-th iteration (we will set $(X^{k+1}, Y^{k+1})$ to be $(U, V)$ if a line search criterion is satisfied; more details can be found in Algorithm \ref{alg_NAUM}). For notational simplicity, we also define
\begin{eqnarray*}
\mathcal{H}_{\alpha}(X,Y,Z):=\frac{\alpha}{2}\|XY^{\top}-Z\|_F^2
\end{eqnarray*}
for any $(X,Y,Z)$. Then, at the $k$-th iteration, we first compute $Z^k$ by \eqref{Zupdate} and, in the line search loop, we compute $U$ in one of the following 3 ways for a given $\mu_k>0$:

\stepcounter{equation}
\begin{itemize}
  \item \textbf{Proximal}
        \begin{equation}\label{Xupdate}
        U \in \mathop{\mathrm{Argmin}}\limits_{X}~\Psi(X)+\mathcal{H}_{\alpha} (X, Y^k, Z^k) + \frac{\mu_k}{2}\|X-X^{k}\|_F^2.\tag{\arabic{section}.\arabic{equation}a}
        \end{equation}
  \item \textbf{Prox-linear}
        \begin{equation}\label{Xlinupdate}
        U \in \mathop{\mathrm{Argmin}}\limits_{X}~\Psi(X) + \langle \nabla_{X}\mathcal{H}_{\alpha}(X^k, Y^k, Z^k),\,X-X^k\rangle
        + \frac{\mu_k}{2}\|X-X^{k}\|_F^2.\tag{\arabic{section}.\arabic{equation}b}
        \end{equation}
  \item \textbf{Hierarchical-prox} ~If $\Psi$ is column-wise separable, i.e., $\Psi(X)=\sum^r_{i=1}\psi_i(\bm{x}_{i})$ for $X=[\bm{x}_1$, $\cdots,\bm{x}_r]\in\mathbb{R}^{m\times r}$, we can update $U$ column-by-column. Specifically, for $i=1,2,\cdots,r$, compute
        \begin{equation}\label{Xhiupdate}
        \bm{u}_i \in \mathop{\mathrm{Argmin}}\limits_{\bm{x}_i}~\psi_i(\bm{x}_i)+\mathcal{H}_{\alpha} (\bm{u}_{j<i}, \bm{x}_i, \bm{x}^k_{j>i}, Y^k, Z^k) + \frac{\mu_k}{2}\|\bm{x}_{i}-\bm{x}_{i}^{k}\|^2, \tag{\arabic{section}.\arabic{equation}c}
        \end{equation}
        where $\bm{u}_{j<i}$ denotes $(\bm{u}_{1},\cdots,\bm{u}_{i-1})$ and $\bm{x}_{j>i}^k$ denotes $(\bm{x}_{i+1}^k,\cdots,\bm{x}_{r}^k)$.
\end{itemize}
After computing $U$, we compute $V$ in one of the following 3 ways for a given $\sigma_k>0$:

\stepcounter{equation}
\begin{itemize}
  \item \textbf{Proximal}
        \begin{equation}\label{Yupdate}
        V \in \mathop{\mathrm{Argmin}}\limits_{Y}~\Phi(Y)+\mathcal{H}_{\alpha} (U, Y, Z^k) + \frac{\sigma_k}{2}\|Y-Y^{k}\|_F^2.\tag{\arabic{section}.\arabic{equation}a}
        \end{equation}
  \item \textbf{Prox-linear}
        \begin{equation}\label{Ylinupdate}
        V \in \mathop{\mathrm{Argmin}}\limits_{Y}~\Phi(Y) + \langle \nabla_{Y}\mathcal{H}_{\alpha}(U, Y^k, Z^k),\,Y-Y^k\rangle
        + \frac{\sigma_k}{2}\|Y-Y^{k}\|_F^2.\tag{\arabic{section}.\arabic{equation}b}
        \end{equation}
  \item \textbf{Hierarchical-prox} ~If $\Phi$ is column-wise separable, i.e., $\Phi(Y)=\sum^r_{i=1}\phi_i(\bm{y}_{i})$ for $Y=[\bm{y}_1$, $\cdots, \bm{y}_r]\in\mathbb{R}^{n\times r}$, we can update $V$ column-by-column. Specifically, for $i=1,2,\cdots,r$, compute
        \begin{equation}\label{Yhiupdate}
        \bm{v}_{i} \in \mathop{\mathrm{Argmin}}\limits_{\bm{y}_i}~\phi_i(\bm{y}_i)+\mathcal{H}_{\alpha} (U, \bm{v}_{j<i}, \bm{y}_i, \bm{y}^k_{j>i}, Z^k) + \frac{\sigma_k}{2}\|\bm{y}_{i}-\bm{y}_i^{k}\|^2, \tag{\arabic{section}.\arabic{equation}c}
        \end{equation}
        where $\bm{v}_{j<i}$ denotes $(\bm{v}_{1},\cdots,\bm{v}_{i-1})$ and $\bm{y}_{j>i}^k$ denotes $(\bm{y}_{i+1}^k,\cdots,\bm{y}_{r}^k)$.
\end{itemize}
For notational simplicity, we further let
\begin{eqnarray}\label{defrho}
\rho := \left\|\mathcal{I}-{\textstyle\frac{\beta}{\alpha+\beta}}\mathcal{A}^*\mathcal{A}\right\|^2
\end{eqnarray}
and let $\gamma\geq0$ be a nonnegative number satisfying
\begin{eqnarray}\label{positiveness}
(\alpha + \gamma) \,\mathcal{I} + \beta \mathcal{A}^*\mathcal{A} \succeq 0.
\end{eqnarray}

\begin{remark}[\textbf{Comments on ``hierarchical-prox"}]
The hierarchical-prox updating scheme requires the column-wise separability of $\Psi$ or $\Phi$. This is satisfied for many common regularizers, for example, $\|\cdot\|_F^2$, $\|\cdot\|_1$, $\|\cdot\|_p^p~(0 < p < 1)$, and the indicator function of the nonnegativity (or box) constraint.
\end{remark}

\begin{remark}[\textbf{Comments on $\rho$ and $\gamma$}]
Since $\mathcal{A}\mathcal{A}^*=\mathcal{I}_q$, we see that the eigenvalues of $\mathcal{A}^*\mathcal{A}$ are either 0 or 1. Then, the eigenvalues of $\mathcal{I}-\frac{\beta}{\alpha+\beta}\mathcal{A}^*\mathcal{A}$ must be either 1 or $\frac{\alpha}{\alpha+\beta}$, and hence $\rho=\max\left\{1,\,\alpha^2/(\alpha+\beta)^2\right\}$. Similarly, the eigenvalues of $-(\alpha \mathcal{I} + \beta \mathcal{A}^*\mathcal{A})$ are either $-\alpha$ or $-(\alpha+\beta)$. Then, \eqref{positiveness} is satisfied whenever $\gamma\geq \max\{0, \,-\alpha, \,-(\alpha+\beta)\}$.
\end{remark}

Now, we are ready to present NAUM as Algorithm \ref{alg_NAUM}.

\begin{algorithm}[ht]
\caption{NAUM for finding a stationary point of $\mathcal{F}$}\label{alg_NAUM}
\textbf{Input:} $(X^0, Y^0)$, $\alpha$ and $\beta$ such that $\frac{1}{\alpha}+\frac1\beta = 1$, $\rho$ as in \eqref{defrho}, $\gamma\geq 0$ satisfying \eqref{positiveness}, $\tau>1$, $c>0$, $\mu^{\min}>0$, $\sigma^{\max} > \sigma^{\min} >0$, and an integer $N\geq0$. Set $k=0$. \vspace{1mm}\\
\textbf{while} a termination criterion is not met, \textbf{do}
\begin{itemize}[leftmargin=2cm]
\item [\textbf{Step 1}.] Compute $Z^k$ by
                  \begin{eqnarray}\label{Zupdate}
                   Z^k = \left(\mathcal{I} - \frac{\beta}{\alpha+\beta}\mathcal{A}^*\mathcal{A}\right)\left(X^k(Y^k)^{\top}\right) + \frac{\beta}{\alpha+\beta}\mathcal{A}^*(\bm{b}).
                  \end{eqnarray}

\item [\textbf{Step 2}.] Choose $\mu_k^0\geq\mu^{\min}$ and $\sigma_k^0 \in [\sigma^{\min}, \,\sigma^{\max}]$ arbitrarily. Set $\tilde{\mu}_k =\mu_k^0$, $\sigma_k = \sigma_k^0$ and $\mu_k^{\max}=(\alpha+2\gamma\rho)\|Y^k\|^2+c$.
         \begin{itemize}[leftmargin=0.8cm]
         \item[\textbf{(2a)}] Set $\mu_k \leftarrow \min\left\{\tilde{\mu}_k, \,\mu_k^{\max}\right\}$. Compute $U$ by either \eqref{Xupdate}, \eqref{Xlinupdate} or \eqref{Xhiupdate}.

         \item[\textbf{(2b)}] Compute $V$ by either \eqref{Yupdate}, \eqref{Ylinupdate} or \eqref{Yhiupdate}.

         \item[\textbf{(2c)}] If
                 \begin{eqnarray}\label{lscond}
                 \begin{aligned}
                 \mathcal{F}(U, V) - \max\limits_{[k-N]_{+}\leq i\leq k}\mathcal{F}(X^i, Y^i)\leq-\frac{c}{2}\left(\|U-X^{k}\|_F^2+\|V-Y^{k}\|_F^2\right),
                 \end{aligned}
                 \end{eqnarray}
                 then go to \textbf{Step 3}.

         \item[\textbf{(2d)}]  If $\mu_k=\mu_k^{\max}$, set $\sigma_k^{\max}=(\alpha+2\gamma\rho)\|U\|^2+c$, $\sigma_k \leftarrow \min\left\{\tau\sigma_k, \,\sigma_k^{\max}\right\}$ and then, go to step \textbf{(2b)}; otherwise, set $\tilde{\mu}_k \leftarrow \tau\mu_k$ and $\sigma_k \leftarrow \tau\sigma_k$ and then, go to step \textbf{(2a)}.
         \end{itemize}

\item [\textbf{Step 3}.] Set $X^{k+1}\leftarrow U$, $Y^{k+1}\leftarrow V$, $\bar{\mu}_k\leftarrow\mu_k$, $\bar{\sigma}_k\leftarrow\sigma_k$, $k \leftarrow k+1$ and go to \textbf{Step 1}.
\end{itemize}
\textbf{end while} \vspace{1mm} \\
\textbf{Output}: $(X^k, Y^k)$ \vspace{1mm}
\end{algorithm}

In Algorithm \ref{alg_NAUM}, the update for $Z^k$ is given explicitly. This is motivated by the condition on $Z$ at a stationary point of $\Theta_{\alpha,\beta}$; see \eqref{optcz}. In fact, following the same arguments in \eqref{zequa}, we see that \eqref{optcz} always holds at $(X^k, Y^k, Z^k)$ with $Z^k$ given in \eqref{Zupdate} when $\mathcal{A}\mathcal{A}^*=\mathcal{I}_q$ and $\frac{1}{\alpha}+\frac{1}{\beta}=1$. If, in addition, $\alpha \mathcal{I} + \beta \mathcal{A}^*\mathcal{A} \succ 0$ holds, one can show that $Z^k$ is actually the optimal solution to the problem $\min_{Z}\left\{\Theta_{\alpha,\beta} (X^k, Y^k, Z)\right\}$. In this case, our NAUM with $N=0$ in \eqref{lscond} can be viewed as an alternating-minimization-based method (see, for example, \cite{abs2013convergence,xy2013a}) applied to the problem $\min_{X,Y,Z}\left\{\Theta_{\alpha,\beta} (X, Y, Z)\right\}$. However, if $\alpha \mathcal{I} + \beta \mathcal{A}^*\mathcal{A}\nsucceq0$,\footnote{This may happen when $0<\alpha<1$ so that $\beta = \alpha(\alpha-1)^{-1}<0$, or $0<\beta<1$ so that $\alpha = \beta(\beta-1)^{-1}<0$.} then the corresponding $\inf_{Z}\left\{\Theta_{\alpha,\beta} (X^k, Y^k, Z)\right\}=-\infty$ for all $k$, and $Z^k$ is only a stationary point of $Z\mapsto \Theta_{\alpha,\beta} (X^k, Y^k, Z)$. In this case, the function value of $\Theta_{\alpha,\beta}$ may increase after updating $Z$ by \eqref{Zupdate}. Fortunately, as we shall see later in \eqref{Ydiff} and \eqref{Xdiff}, as long as $\mathcal{A}\mathcal{A}^*=\mathcal{I}_q$ and $\frac{1}{\alpha}+\frac{1}{\beta}=1$, we still have $\Theta_{\alpha,\beta} (X^{k+1}, Y^{k+1}, Z^k) < \Theta_{\alpha,\beta} (X^k, Y^k, Z^{k})$ by updating $X^{k+1}$ and $Y^{k+1}$ with properly chosen parameters $\mu_k$ and $\sigma_k$. Thus, if the possible increase in $\Theta_{\alpha,\beta}$ induced by the $Z$-update is not too large, one can still ensure $\Theta_{\alpha,\beta} (X^{k+1}, Y^{k+1}, Z^{k+1}) < \Theta_{\alpha,\beta} (X^k, Y^k, Z^{k})$. Moreover, it can be seen from Lemma \ref{supplem} and \eqref{Zupdate} that $\mathcal{F}(X^k, Y^k) = \Theta_{\alpha,\beta}(X^k,Y^k$, $Z^k)$ and hence the decrease of $\Theta_{\alpha,\beta}$ translates to that of $\cal F$ (see Lemma \ref{suffdes} below). In view of this, $\Theta_{\alpha,\beta}$ is a valid potential function for minimizing $\cal F$ as long as $\mathcal{A}\mathcal{A}^*=\mathcal{I}_q$ and $\frac{1}{\alpha}+\frac{1}{\beta}=1$, even when $\beta<0$ or $\alpha<0$. Allowing negative $\alpha$ or $\beta$ makes our NAUM (even with $N=0$ in \eqref{lscond}) different from the classical alternating minimization schemes.

Our NAUM also allows $U$ and $V$ to be updated in three different ways respectively, and hence there are 9 possible combinations. Thus, one can choose suitable updating schemes to fit different applications. In particular, if $\Psi$ or $\Phi$ are column-wise separable, taking advantage of the structure of $\Theta_{\alpha,\beta}$ and the fact that $XY^{\top}$ can be written as $\sum^r_{i=1}\bm{x}_i\bm{y}_i^{\top}$ with $X=[\bm{x}_1, \cdots, \bm{x}_r]\in\mathbb{R}^{m\times r}$ and $Y=[\bm{y}_1, \cdots, \bm{y}_r]\in\mathbb{R}^{n\times r}$, one can update $X$ or $Y$ column-wise even when $\mathcal{A}\neq\mathcal{I}$. The motivation for updating $X$ (or $Y$) column-wise rather than updating the whole $X$ (or $Y$) is that the resulting subproblems \eqref{Xhiupdate} (or \eqref{Yhiupdate}) can be reduced to the computation of the proximal mapping of $\psi_i$ (or $\phi_i$), which is easy for many commonly used $\psi_i$ (or $\phi_i$). Indeed, from \eqref{Xhiupdate} and \eqref{Yhiupdate}, $\bm{u}_i$ and $\bm{v}_i$ are given by
\begin{eqnarray}\label{subprouv}
\left\{\begin{aligned}
&\bm{u}_i \in \mathop{\mathrm{Argmin}}\limits_{\bm{x}_i}\,\left\{\psi_i(\bm{x}_{i}) +\frac{\alpha}{2}\left\|\bm{x}_i(\bm{y}^k_i)^{\top}-P_i^k\right\|_F^2 + \frac{\mu_k}{2}\|\bm{x}_i - \bm{x}_i^k\|^2\right\},   \\
&\bm{v}_i \in \mathop{\mathrm{Argmin}}\limits_{\bm{y}_i}\,\left\{\phi_i(\bm{y}_{i}) +\frac{\alpha}{2}\left\|\bm{u}_i\bm{y}_i^{\top}-Q_i^k\right\|_F^2 + \frac{\sigma_k}{2}\|\bm{y}_i - \bm{y}_i^k\|^2\right\},
\end{aligned}\right.
\end{eqnarray}
where $P_i^k$ and $Q_i^k$ are defined by
\begin{eqnarray}\label{defPQ}
\begin{aligned}
&P_i^k:=~Z^k-{\textstyle\sum^{i-1}_{j=1}}\bm{u}_j(\bm{y}^k_j)^{\top}
-{\textstyle\sum^r_{j=i+1}}\bm{x}^{k}_j(\bm{y}^k_j)^{\top}, \\
&Q_i^k:=~Z^k-{\textstyle\sum^{i-1}_{j=1}}\bm{u}_j\bm{v}_j^{\top}
-{\textstyle\sum^r_{j=i+1}}\bm{u}_j(\bm{y}^k_j)^{\top}.
\end{aligned}
\end{eqnarray}
Then, from Proposition \ref{propuse}, we can reformulate the subproblems in \eqref{subprouv} and obtain the corresponding solutions by computing the proximal mappings of $\psi_i$ and $\phi_i$, which can be computed efficiently when $\psi_i$ and $\phi_i$ are some common regularizers used in the literature. In particular, when $\psi_i(\cdot)$ and $\phi_i(\cdot)$ are $\|\cdot\|_1$, $\|\cdot\|_2^2$ or the indicator function of the box constraint, these subproblems have closed-form solutions. This updating strategy has also been used for NMF; see, for example, \cite{cza2007hierarchical,lz2009fastnmf,llwy2012sparse}. However, the methods used in \cite{cza2007hierarchical,lz2009fastnmf,llwy2012sparse} can only be applied for some specific problems with $\mathcal{A}=\mathcal{I}$, while NAUM can be applied for more general problems with $\mathcal{A}\mathcal{A}^*=\mathcal{I}_q$.

Our NAUM adapts a non-monotone line search criterion (see Step 2 in Algorithm \ref{alg_NAUM}) to improve the numerical performance. This is motivated by recent studies on non-monotone algorithms with promising performances; see, for example, \cite{clp2016penalty,gzlhy2013a,wnf2009sparse}. However, different from the non-monotone line search criteria used there, NAUM only includes $(U, V)$ in the line search loop and checks the stopping criterion \eqref{lscond} after updating a pair of $(U, V)$, rather than checking \eqref{lscond} immediately once $U$ or $V$ is updated. Thus, we do not need to compute the function value after updating each block of variable. This may reduce the cost of the line search and make NAUM more practical, especially when computing the function value is relatively expensive.

Before moving to the convergence analysis of NAUM, we would like to point out an interesting connection between NAUM and the low-rank matrix fitting algorithm, LMaFit \cite{wyz2012solving}, for solving the following matrix completion model without regularizers:
\begin{eqnarray*}
\min\limits_{X,Y}~\frac{1}{2}\left\|\mathcal{P}_{\Omega}(XY^{\top}-M)\right\|_F^2,
\end{eqnarray*}
where $\Omega$ is the index set of the known entries of $M$, and $\mathcal{P}_{\Omega}(Z)$ keeps the entries of $Z$ in $\Omega$ and sets the remaining ones to zero. If we apply our NAUM with \eqref{Xupdate} and \eqref{Yupdate}, then at the $k$-th iteration, the iterates $Z^k$, $X^{k+1}$ and $Y^{k+1}$ are given by
\begin{eqnarray*}
\begin{aligned}
Z^k &= \left(\mathcal{I}-{\textstyle\frac{\beta}{\alpha+\beta}}\mathcal{P}_{\Omega}\right)X^k(Y^k)^{\top} + {\textstyle\frac{\beta}{\alpha+\beta}}\mathcal{P}_{\Omega}(M),     \\
X^{k+1} &= \left(\bar{\mu}_kX^k + \alpha Z^kY^k\right)\left(\bar{\mu}_kI+\alpha(Y^k)^{\top}Y^k\right)^{-1},             \\
Y^{k+1} &= \left(\bar{\sigma}_kY^k + \alpha (Z^k)^{\top}X^{k+1}\right)\left(\bar{\sigma}_kI+\alpha(X^{k+1})^{\top}X^{k+1}\right)^{-1}.   \\
\end{aligned}
\end{eqnarray*}
One can verify that the sequence $\{(Z^k, X^{k+1}, Y^{k+1})\}$ above can be equivalently generated by the following scheme with
$\widetilde{Z}^{0}=\mathcal{P}_{\Omega}(M) + \mathcal{P}_{\Omega^c}\left(X^{0}(Y^{0})^{\top}\right)$:
\begin{eqnarray*}
\begin{aligned}
Z^k &= {\textstyle\frac{\beta}{\alpha+\beta}}\widetilde{Z}^{k}
+\left(1-{\textstyle\frac{\beta}{\alpha+\beta}}\right)X^k(Y^k)^{\top} ,     \\
X^{k+1} &= \left(\bar{\mu}_kX^k + \alpha Z^kY^k\right)\left(\bar{\mu}_kI+\alpha(Y^k)^{\top}Y^k\right)^{-1},             \\
Y^{k+1} &= \left(\bar{\sigma}_kY^k + \alpha (Z^k)^{\top}X^{k+1}\right)\left(\bar{\sigma}_kI+\alpha(X^{k+1})^{\top}X^{k+1}\right)^{-1},    \\
\widetilde{Z}^{k+1} &= \mathcal{P}_{\Omega}(M) + \mathcal{P}_{\Omega^c}\left(X^{k+1}(Y^{k+1})^{\top}\right), \end{aligned}
\end{eqnarray*}
where $\Omega^c$ is the complement set of $\Omega$. Surprisingly, when $\bar{\mu}_k=\bar{\sigma}_k=0$, this scheme is exactly the SOR(successive over-relaxation)-like scheme used in LMaFit (see \cite[Eq.(2.11)]{wyz2012solving}) with $\omega:=\frac{\beta}{\alpha+\beta}$ being an over-relaxation weight. With this connection, our NAUM, in some sense, can be viewed as an SOR-based algorithm. Moreover, just like the classical SOR for solving
a system of linear equations, LMaFit with $\omega>1$ also appears to be more efficient from the extensive numerical experiments reported in \cite{wyz2012solving}. Then, it is natural to consider $\frac{\beta}{\alpha+\beta}>1$ and hence $\frac{1}{\alpha}>1$ (since $\frac{1}{\alpha}+\frac{1}{\beta}=1$) in NAUM. This also gives some insights for the necessity of allowing more flexibilities in choosing $\alpha$ and $\beta$, and the promising performance of NAUM with a relatively small $\alpha\in(0,1)$ as we shall see in Section \ref{secnum}.

\section{Convergence analysis of NAUM}\label{secconv}

In this section, we discuss the convergence properties of Algorithm \ref{alg_NAUM}. First, we present the first-order optimality conditions for the three different updating schemes in (2a) of Algorithm \ref{alg_NAUM} as follows:

\stepcounter{equation}
\begin{itemize}
\item \textbf{Proximal}
      \begin{equation}\label{Xoptcon}
      0 \in \partial\Psi(U) + \alpha\left(U(Y^{k})^{\top} - Z^k\right)Y^k + \mu_k(U - X^k). \tag{\arabic{section}.\arabic{equation}a}
      \end{equation}

\item \textbf{Prox-linear}
      \begin{equation}\label{Xlinoptcon}
      0 \in \partial\Psi(U) + \alpha\left(X^k(Y^{k})^{\top} - Z^k\right)Y^k + \mu_k(U - X^k). \tag{\arabic{section}.\arabic{equation}b}
      \end{equation}

\item \textbf{Hierarchical-prox} ~For $i = 1,2,\cdots,r$,
      \begin{equation}\label{Xhioptcon}
      0 \in \partial\psi_i(\bm{u}_{i}) +\alpha\left({\textstyle\sum^i_{j=1}}\bm{u}_j(\bm{y}^k_j)^{\top}
      +{\textstyle\sum^r_{j=i+1}}\bm{x}^{k}_j(\bm{y}^k_j)^{\top} - Z^k\right)\bm{y}_i^k + \mu_k(\bm{u}_i - \bm{x}_i^k). \tag{\arabic{section}.\arabic{equation}c}
      \end{equation}
\end{itemize}
Similarly, the first-order optimality conditions for the three different updating schemes in (2b) of Algorithm \ref{alg_NAUM} are

\stepcounter{equation}
\begin{itemize}
\item \textbf{Proximal}
      \begin{equation}\label{Yoptcon}
      0 \in \partial\Phi(V)+\alpha\left(UV^{\top} - Z^k\right)^{\top}U + \sigma_k(V - Y^k).
      \tag{\arabic{section}.\arabic{equation}a}
      \end{equation}

\item \textbf{Prox-linear}
      \begin{equation}\label{Ylinoptcon}
      0 \in \partial\Phi(V)+\alpha\left(U(Y^k)^{\top} - Z^k\right)^{\top}U + \sigma_k(V - Y^k).
      \tag{\arabic{section}.\arabic{equation}b}
      \end{equation}

\item \textbf{Hierarchical-prox} ~For $i = 1,2,\cdots,r$,
      \begin{equation}\label{Yhioptcon}
      0 \in \partial\phi_i(\bm{v}_{i})+\alpha\left({\textstyle\sum^i_{j=1}}\bm{u}_j\bm{v}_j^{\top}
      +{\textstyle\sum^r_{j=i+1}}\bm{u}_j(\bm{y}^k_j)^{\top} - Z^k\right)^{\top}\bm{u}_i + \sigma_k(\bm{v}_i - \bm{y}_i^k). \tag{\arabic{section}.\arabic{equation}c}
      \end{equation}
\end{itemize}

We also need to make the following assumptions.

\begin{assumption}\label{assumfun}
\indent
\begin{itemize}[leftmargin=0.9cm]
\item[{\bf(a1)}] $\Psi$, $\Phi$ are proper, closed, level-bounded functions and continuous on their domains respectively;

\item[{\bf(a2)}] $\mathcal{A}\mathcal{A}^* = \mathcal{I}_q$;

\item[{\bf(a3)}] $\frac{1}{\alpha} + \frac{1}{\beta} = 1$.
\end{itemize}
\end{assumption}

\begin{remark}
(i) From {\rm(a1)}, one can see from \cite[Theorem~1.9]{rw1998variational} that $\inf\Psi$ and $\inf\Phi$ are finite, i.e., $\Psi$ and $\Phi$ are bounded from below. In particular, the iterates \eqref{Xupdate}, \eqref{Xlinupdate}, \eqref{Xhiupdate}, \eqref{Yupdate}, \eqref{Ylinupdate} and \eqref{Yhiupdate} are well defined; (ii) The continuity assumption in {\rm (a1)} holds for many common regularizers, for example, $\ell_1$-norm, nuclear norm and the indicator function of a nonempty closed set; (iii) {\rm (a2)} is satisfied for some commonly used linear maps, for example, the identity map and the sampling map.
\end{remark}

We start our convergence analysis by proving the following auxiliary lemma.

\begin{lemma}[\textbf{Sufficient descent of $\mathcal{F}$}]\label{suffdes}
Suppose that Assumption \ref{assumfun} holds. Let $(X^k, Y^k)$ be generated by Algorithm \ref{alg_NAUM} at the $k$-th iteration, and $(U, V)$ be the candidate for $(X^{k+1}, Y^{k+1})$ generated by steps {\bf(2a)} and {\bf(2b)}. Then, for any integer $k \geq 0$, we have
\begin{eqnarray}\label{succhange}
\mathcal{F}(U, V) - \mathcal{F}(X^k, Y^k)
\leq -\frac{\mu_k\!-\!(\alpha\!+\!2\gamma\rho)\|Y^{k}\|^2}{2}\,\|U-X^{k}\|_F^2
-\frac{\sigma_k\!-\!(\alpha\!+\!2\gamma\rho)\|U\|^2}{2}\,\|V-Y^{k}\|_F^2.
\end{eqnarray}
\end{lemma}
\beginproof
First, from Lemma \ref{supplem} and \eqref{Zupdate}, we see that $\mathcal{F}(X^k, Y^k) = \Theta_{\alpha,\beta}(X^k,Y^k,Z^k)$. For any $(U,V)$, let
\begin{eqnarray}\label{wupdate}
W = {\textstyle\left(\mathcal{I} - \frac{\beta}{\alpha+\beta}\mathcal{A}^*\mathcal{A}\right)}\left(UV^{\top}\right) + {\textstyle\frac{\beta}{\alpha+\beta}}\mathcal{A}^*(\bm{b}).
\end{eqnarray}
Then, from Lemma \ref{supplem}, we have $\mathcal{F}(U, V) = \Theta_{\alpha,\beta}(U,V,W)$. Thus, to establish \eqref{succhange}, we only need to consider the difference $\Theta_{\alpha,\beta}(U,V,W)-\Theta_{\alpha,\beta}(X^k,Y^k,Z^k)$.

We start by noting that
\begin{eqnarray}\label{AtAZ}
\begin{aligned}
\mathcal{A}^*\mathcal{A}(W)
&=\left(\mathcal{A}^*\mathcal{A} - {\textstyle\frac{\beta}{\alpha+\beta}}\mathcal{A}^*\left(\mathcal{A}\mathcal{A}^*\right)\mathcal{A}\right)\left(UV^{\top}\right)+ {\textstyle\frac{\beta}{\alpha+\beta}}\mathcal{A}^*\left(\mathcal{A}\mathcal{A}^*\right)(\bm{b}) \\
&={\textstyle\frac{\alpha}{\alpha+\beta}}\mathcal{A}^*\mathcal{A}\left(UV^{\top}\right)+{\textstyle\frac{\beta}{\alpha+\beta}}\mathcal{A}^*(\bm{b}),
\end{aligned}
\end{eqnarray}
where the last equality follows from (a2) in Assumption \ref{assumfun}. Then, we obtain that
\begin{eqnarray*}
\begin{aligned}
&~~\nabla_{Z} \Theta_{\alpha,\beta} (U, V, W)
=\alpha(W-UV^{\top})+\beta\mathcal{A}^*\mathcal{A}(W)-\beta\mathcal{A}^*(\bm{b}) \\
&= \alpha\left[-{\textstyle\frac{\beta}{\alpha+\beta}}\mathcal{A}^*\mathcal{A}(UV^{\top})
+{\textstyle\frac{\beta}{\alpha+\beta}}\mathcal{A}^*(\bm{b})\right]
+\beta\left[{\textstyle\frac{\alpha}{\alpha+\beta}}\mathcal{A}^*\mathcal{A}\left(UV^{\top}\right)+
{\textstyle\frac{\beta}{\alpha+\beta}}\mathcal{A}^*(\bm{b})\right] - \beta\mathcal{A}^*(\bm{b})
= 0,
\end{aligned}
\end{eqnarray*}
where the second equality follows from \eqref{wupdate} and \eqref{AtAZ}. Moreover, since $\gamma$ is chosen such that $(\alpha + \gamma) \mathcal{I} + \beta \mathcal{A}^*\mathcal{A} \succeq 0$ (see \eqref{positiveness}), we see that, for any $k\geq0$, the function $Z \longmapsto \Theta_{\alpha,\beta} (U, V, Z)+\frac{\gamma}{2}\|Z - Z^{k}\|_F^2$ is convex and hence
\begin{eqnarray*}
\begin{aligned}
&\quad \Theta_{\alpha,\beta} (U, V, Z^{k})+\underbrace{\frac{\gamma}{2}\|Z^{k}-Z^{k}\|_F^2}_{=0} \\
&\geq \Theta_{\alpha,\beta} (U, V, W)+\frac{\gamma}{2}\|W-Z^{k}\|_F^2
+\langle\underbrace{\nabla_{Z} \Theta_{\alpha,\beta} (U, V, W)}_{=0}\,+\,\gamma(W-Z^k), \,Z^k-W\rangle,
\end{aligned}
\end{eqnarray*}
which implies that
\begin{eqnarray}\label{addineq1}
\Theta_{\alpha,\beta}(U, V, W) - \Theta_{\alpha,\beta}(U, V, Z^{k}) \leq \frac{\gamma}{2}\|W - Z^{k}\|_F^2.
\end{eqnarray}
Then, substituting \eqref{Zupdate} and \eqref{wupdate} into \eqref{addineq1}, we obtain
\begin{eqnarray}\label{Zdiff}
\begin{aligned}
&~~\Theta_{\alpha,\beta}(U, V, W) - \Theta_{\alpha,\beta}(U, V, Z^{k}) \\
&\leq\frac{\gamma}{2}\left\|{\textstyle\left(\mathcal{I}-\frac{\beta}{\alpha+\beta}\mathcal{A}^*\mathcal{A}\right)}
\left(UV^{\top}-X^{k}(Y^{k})^{\top}\right)\right\|_F^2
~\leq~ \frac{\gamma}{2}\left\|{\textstyle\mathcal{I} - \frac{\beta}{\alpha+\beta}\mathcal{A}^*\mathcal{A}}\right\|^2
\cdot\left\|UV^{\top}-X^{k}(Y^{k})^{\top}\right\|_F^2  \\
&=\frac{\gamma\rho}{2}\left\|U(V-Y^k)^{\top}+(U-X^{k})(Y^{k})^{\top}\right\|_F^2
~\leq~\frac{\gamma\rho}{2}\Big{(}\left\|U(V-Y^k)^{\top}\right\|_F
+\left\|(U-X^{k})(Y^{k})^{\top}\right\|_F\Big{)}^2 \\
&\stackrel{(\mathrm{i})}{\leq} \frac{\gamma\rho}{2}\Big{(}\|U\|\|V-Y^k\|_F+\|Y^{k}\|\|U-X^{k}\|_F\Big{)}^2
~\stackrel{(\mathrm{ii})}{\leq}~\gamma\rho\Big{(}\|U\|^2\|V-Y^k\|^2_F+\|Y^{k}\|^2\|U-X^{k}\|^2_F\Big{)},
\end{aligned}
\end{eqnarray}
where the equality follows from the definition of $\rho$ in \eqref{defrho}; (i) follows from the relation $\|AB\|_F \leq \|A\|\|B\|_F$; and (ii) follows from the relation $\|a+b\|^2\leq2\|a\|^2+2\|b\|^2$.

Next, we show that
\begin{eqnarray}\label{Ydiff}
\Theta_{\alpha,\beta}(U, V, Z^{k}) - \Theta_{\alpha,\beta}(U, Y^k, Z^k)\leq \frac{\alpha\|U\|^2-\sigma_k}{2}\|V-Y^{k}\|_F^2.
\end{eqnarray}
To this end, we consider the following three cases.
\begin{itemize}
  \item Proximal: In this case, we have
      \begin{eqnarray*}
      \begin{aligned}
      &~~\Theta_{\alpha,\beta}(U, V, Z^{k}) - \Theta_{\alpha,\beta}(U, Y^k, Z^k)
      =\Phi(V)+\mathcal{H}_{\alpha}(U, V, Z^{k}) - \Phi(Y^k) - \mathcal{H}_{\alpha}(U, Y^k, Z^{k})\\
      &=\left[\Phi(V)+\mathcal{H}_{\alpha}(U, V, Z^{k})+\frac{\sigma_k}{2}\|V-Y^{k}\|_F^2\right]
       -\left[\Phi(Y^k)+\mathcal{H}_{\alpha}(U, Y^k, Z^k)\right]-\frac{\sigma_k}{2}\|V-Y^{k}\|_F^2\\
      &\leq -\frac{\sigma_k}{2}\|V-Y^{k}\|_F^2,
      \end{aligned}
      \end{eqnarray*}
      where the inequality follows from the definition of $V$ as a minimizer of \eqref{Yupdate}. This implies \eqref{Ydiff}.

  \item Prox-linear: In this case, we have
      \begin{eqnarray*}
      \begin{aligned}
      &~~\Theta_{\alpha,\beta}(U, V, Z^{k}) - \Theta_{\alpha,\beta}(U, Y^k, Z^k)
      =\Phi(V)+\mathcal{H}_{\alpha}(U, V, Z^{k}) - \Phi(Y^k) - \mathcal{H}_{\alpha}(U, Y^k, Z^{k}) \\
      &\leq \Phi(V) + \mathcal{H}_{\alpha}(U, Y^k, Z^{k})+ \langle \nabla_{Y}\mathcal{H}_{\alpha}(U, Y^k, Z^k),\,V-Y^k\rangle+\frac{\alpha\|U\|^2}{2}\|V-Y^{k}\|_F^2 \\
      &\quad -\Phi(Y^k)-\mathcal{H}_{\alpha}(U, Y^k, Z^{k}) \\
      &=\Phi(V)+ \langle \nabla_{Y}\mathcal{H}_{\alpha}(U, Y^k, Z^k),\,V-Y^k\rangle+\frac{\sigma_k}{2}\|V-Y^{k}\|_F^2-\Phi(Y^k) + \frac{\alpha\|U\|^2-\sigma_k}{2}\|V-Y^{k}\|_F^2 \\
      &\leq \frac{\alpha\|U\|^2-\sigma_k}{2}\|V-Y^{k}\|_F^2,
      \end{aligned}
      \end{eqnarray*}
      where the first inequality follows from the fact that $Y\mapsto\nabla_{Y}\mathcal{H}_{\alpha}(X, Y, Z)$ is Lipschitz with modulus $\alpha\|X\|^2$ and the last inequality follows from the definition of $V$ as a minimizer of \eqref{Ylinupdate}.

  \item Hierarchical-prox: In this case, for any $1\leq i \leq r$, we have
        \begin{eqnarray*}
        \begin{aligned}
        &~~\Theta_{\alpha,\beta} (U, \bm{v}_{j<i}, \bm{v}_i, \bm{y}^k_{j>i}, Z^k)
        -\Theta_{\alpha,\beta} (U, \bm{v}_{j<i}, \bm{y}^{k}_i, \bm{y}^k_{j>i}, Z^k) \\
        &=\phi_i(\bm{v}_i)+\mathcal{H}_{\alpha} (U, \bm{v}_{j<i}, \bm{v}_i, \bm{y}^k_{j>i}, Z^k)-
        \phi_i(\bm{y}_i^k)-\mathcal{H}_{\alpha} (U, \bm{v}_{j<i}, \bm{y}^{k}_i, \bm{y}^k_{j>i}, Z^k)  \\
        &=\left[\phi_i(\bm{v}_i)+\mathcal{H}_{\alpha} (U, \bm{v}_{j<i}, \bm{v}_i, \bm{y}^k_{j>i}, Z^k) + \frac{\sigma_k}{2}\|\bm{v}_{i}-\bm{y}^{k}_i\|^2\right] - \frac{\sigma_k}{2}\|\bm{v}_{i}-\bm{y}^{k}_i\|^2 \\
        &\quad-\left[\phi_i(\bm{y}_i^k)+\mathcal{H}_{\alpha} (U, \bm{v}_{j<i}, \bm{y}^{k}_i, \bm{y}^k_{j>i}, Z^k)\right] \\
        &\leq -\frac{\sigma_k}{2}\|\bm{v}_{i}-\bm{y}^{k}_i\|^2,
        \end{aligned}
        \end{eqnarray*}
        where the inequality follows from the definition of $\bm{v}_i$ as a minimizer of
        \eqref{Yhiupdate}. Then, summing the above relation from $i=r$ to $i=1$ and simplifying the resulting inequality, we obtain \eqref{Ydiff}.
\end{itemize}

Similarly, we can show that
\begin{eqnarray}\label{Xdiff}
\Theta_{\alpha,\beta}(U, Y^{k}, Z^{k}) - \Theta_{\alpha,\beta}(X^{k}, Y^k, Z^k)
\leq \frac{\alpha\|Y^k\|^2-\mu_k}{2}\|U-X^{k}\|_F^2
\end{eqnarray}
by considering the following three cases.
\begin{itemize}
  \item Proximal: In this case, we have
      \begin{eqnarray*}
      \begin{aligned}
      &~~\Theta_{\alpha,\beta}(U, Y^{k}, Z^{k}) - \Theta_{\alpha,\beta}(X^{k}, Y^k, Z^k)
      =\Psi(U)+\mathcal{H}_{\alpha}(U, Y^{k}, Z^{k}) - \Psi(X^k) - \mathcal{H}_{\alpha}(X^{k}, Y^k, Z^{k})\\
      &=\left[\Psi(U)+\mathcal{H}_{\alpha}(U, Y^{k}, Z^{k})+\frac{\mu_k}{2}\|U-X^{k}\|_F^2\right]-\left[\Psi(X^k)+\mathcal{H}_{\alpha}(X^{k}, Y^k, Z^k)\right]-\frac{\mu_k}{2}\|U-X^{k}\|_F^2 \\
      &\leq -\frac{\mu_k}{2}\|U-X^{k}\|_F^2,
      \end{aligned}
      \end{eqnarray*}
      where the inequality follows from the definition of $U$ as a minimizer of \eqref{Xupdate}. This implies \eqref{Xdiff}.

  \item Prox-linear: In this case, we have
      \begin{eqnarray*}
      \begin{aligned}
      &~~\Theta_{\alpha,\beta}(U, Y^{k}, Z^{k}) - \Theta_{\alpha,\beta}(X^{k}, Y^k, Z^k)
      =\Psi(U)+\mathcal{H}_{\alpha}(U, Y^{k}, Z^{k}) - \Psi(X^k) - \mathcal{H}_{\alpha}(X^{k}, Y^k, Z^{k}) \\
      &\leq \Psi(U) + \mathcal{H}_{\alpha}(X^k, Y^k, Z^{k})+ \langle \nabla_{X}\mathcal{H}_{\alpha}(X^k, Y^k, Z^k),\,U-X^k\rangle+\frac{\alpha\|Y^k\|^2}{2}\|U-X^{k}\|_F^2  \\
      &\quad -\Psi(X^k)-\mathcal{H}_{\alpha}(X^k, Y^k, Z^{k}) \\
      &=\Psi(U)+ \langle \nabla_{X}\mathcal{H}_{\alpha}(X^k, Y^k, Z^k),\,U-X^k\rangle+\frac{\mu_k}{2}\|U-X^{k}\|_F^2-\Psi(X^k) + \frac{\alpha\|Y^k\|^2-\mu_k}{2}\|U-X^{k}\|_F^2 \\
      &\leq \frac{\alpha\|Y^k\|^2-\mu_k}{2}\|U-X^{k}\|_F^2,
      \end{aligned}
      \end{eqnarray*}
      where the first inequality follows from the fact that $\nabla_{X}\mathcal{H}_{\alpha}(X, Y, Z)$ is Lipschitz with modulus $\alpha\|Y\|^2$ and the last inequality follows from the definition of $U$ as a minimizer of \eqref{Xlinupdate}.

  \item Hierarchical-prox: In this case, for any $1\leq i \leq r$, we have
        \begin{eqnarray*}
        \begin{aligned}
        &~~\Theta_{\alpha,\beta} (\bm{u}_{j<i}, \bm{u}_i, \bm{x}^k_{j>i}, Y^k, Z^k)
        -\Theta_{\alpha,\beta} (\bm{u}_{j<i}, \bm{x}^{k}_i, \bm{x}^k_{j>i}, Y^k, Z^k) \\
        &=\psi_i(\bm{u}_i)+\mathcal{H}_{\alpha} (\bm{u}_{j<i}, \bm{u}_i, \bm{x}^k_{j>i}, Y^k, Z^k)
        -\psi_i(\bm{x}_i^k)-\mathcal{H}_{\alpha} (\bm{u}_{j<i}, \bm{x}^{k}_i, \bm{x}^k_{j>i}, Y^k, Z^k) \\
        &=\left[\psi_i(\bm{u}_i)+\mathcal{H}_{\alpha} (\bm{u}_{j<i}, \bm{u}_i, \bm{x}^k_{j>i}, Y^k, Z^k)
        + \frac{\mu_k}{2}\|\bm{u}_{i}-\bm{x}^{k}_i\|^2\right] - \frac{\mu_k}{2}\|\bm{u}_{i}-\bm{x}^{k}_i\|^2 \\
        &\quad-\left[\psi_i(\bm{x}_i^k)+\mathcal{H}_{\alpha} (\bm{u}_{j<i}, \bm{x}^{k}_i, \bm{x}^k_{j>i}, Y^k, Z^k) \right] \\
        &\leq -\frac{\mu_k}{2}\|\bm{u}_{i}-\bm{x}^{k}_i\|^2,
        \end{aligned}
        \end{eqnarray*}
        where the inequality follows from the definition of $\bm{u}_i$ as a minimizer of \eqref{Xhiupdate}. Then, summing the above relation from $i=r$ to $i=1$ and simplifying the resulting inequality, we obtain \eqref{Xdiff}.

\end{itemize}

Now, summing \eqref{Zdiff}, \eqref{Ydiff} and \eqref{Xdiff}, and using $\mathcal{F}(U,V) = \Theta_{\alpha,\beta}(U,V,W)$ and $\mathcal{F}(X^k, Y^k) =\Theta_{\alpha,\beta}(X^k$, $Y^k,Z^k)$, we obtain \eqref{succhange}. This completes the proof.
\endproof

From the above lemma, we see that the sufficient descent of $\mathcal{F}(X, Y)$ can be guaranteed as long as $\mu_k$ and $\sigma_k$ are sufficiently large. Thus, based on this lemma, we can show in the following proposition that our non-monotone line search criterion \eqref{lscond} in Algorithm \ref{alg_NAUM} is well defined.

\begin{proposition}[\textbf{Well-definedness of the line search criterion}]\label{guals}
Suppose that Assumption \ref{assumfun} holds and Algorithm~\ref{alg_NAUM} is applied. Then, for each $k \geq 0$, the line search criterion \eqref{lscond} is satisfied after finitely many inner iterations.
\end{proposition}
\beginproof
We prove this proposition by contradiction. Assume that there exists a $k\geq 0$ such that the line search criterion \eqref{lscond} cannot be satisfied after finitely many inner iterations. Note from (2a) and (2d) in Step 2 of Algorithm \ref{alg_NAUM} that $\mu_k \leq \mu_k^{\max}=(\alpha+2\gamma\rho)\|Y^k\|^2+c$ and hence $\mu_k = \mu_k^{\max}$ must be satisfied after finitely many inner iterations. Let $n_k$ denote the number of inner iterations when $\mu_k = \mu_k^{\max}$ is satisfied for the \textit{first} time. If $\mu^0_k \geq \mu_k^{\max}$, then $n_k=1$; otherwise, we have
\begin{eqnarray*}
\mu^{\min} \tau^{n_k-2} \leq \mu^0_k \tau^{n_k-2} < \mu_k^{\max},
\end{eqnarray*}
which implies that
\begin{eqnarray}\label{upbnk}
n_k \leq \left\lfloor\frac{\log(\mu_k^{\max})-\log(\mu^{\min})}{\log\tau} + 2 \right\rfloor.
\end{eqnarray}
Then, from (2d) in Step 2 of Algorithm \ref{alg_NAUM}, we have $U\equiv U_{\mu^{\max}_k}$ and $\sigma^{\max}_k=(\alpha+2\gamma\rho)\|U_{\mu^{\max}_k}\|^2+c$ after at most $n_k+1$ inner iterations, where $U_{\mu^{\max}_k}$ is computed by \eqref{Xupdate}, \eqref{Xlinupdate} or \eqref{Xhiupdate} with $\mu_k = \mu_k^{\max}$. Moreover, we see that $\sigma_k=\sigma^{\max}_k$ must be satisfied after finitely many inner iterations. Similarly, let $\hat{n}_k$ denote the number of inner iterations when $\sigma_k=\sigma^{\max}_k$ is satisfied for the \textit{first} time. If $\sigma^0_k > \sigma_k^{\max}$, then $\hat{n}_k=n_k$; if $\sigma^0_k = \sigma_k^{\max}$, then $\hat{n}_k=0$; otherwise, we have
\begin{eqnarray*}
\sigma^{\min} \tau^{\hat{n}_k-1} \leq \sigma^0_k \tau^{\hat{n}_k-1} < \sigma_k^{\max},
\end{eqnarray*}
which implies that
\begin{eqnarray*}
\hat{n}_k \leq \left\lfloor\frac{\log(\sigma_k^{\max})-\log(\sigma^{\min})}{\log\tau} + 1 \right\rfloor.
\end{eqnarray*}
Thus, after at most $\max\{n_k, \,\hat{n}_k\}+1$ inner iterations, we must have $V\equiv V_{\sigma^{\max}_k}$, where $V_{\sigma^{\max}_k}$ is computed by \eqref{Yupdate}, \eqref{Ylinupdate} or \eqref{Yhiupdate} with $\sigma_k = \sigma_k^{\max}$. Therefore, after at most $\max\{n_k, \,\hat{n}_k\}+1$ inner iterations, we have
\begin{eqnarray*}
\begin{aligned}
&~~\mathcal{F}(U_{\mu^{\max}_k}, V_{\sigma^{\max}_k}) - \mathcal{F}(X^k, Y^k)   \\
&\leq-\frac{\mu^{\max}_k-(\alpha+2\gamma\rho)\|Y^{k}\|^2}{2}\,\|U_{\mu^{\max}_k}-X^{k}\|_F^2
-\frac{\sigma^{\max}_k-(\alpha+2\gamma\rho)\|U_{\mu^{\max}_k}\|^2}{2}\,\|V_{\sigma^{\max}_k}-Y^{k}\|_F^2 \\
&=-\frac{c}{2}\left(\|U_{\mu^{\max}_k}-X^{k}\|_F^2+\|V_{\sigma^{\max}_k}-Y^{k}\|_F^2\right),
\end{aligned}
\end{eqnarray*}
where the inequality follows from \eqref{succhange} and the equality follows from $\mu_k^{\max}=(\alpha+2\gamma\rho)\|Y^k\|^2+c$ and $\sigma^{\max}_k=(\alpha+2\gamma\rho)\|U_{\mu^{\max}_k}\|^2+c$. This together with
\begin{eqnarray*}
\mathcal{F}(X^k, Y^k) \leq \max\limits_{[k-N]_{+}\leq i\leq k}\mathcal{F}(X^i, Y^i)
\end{eqnarray*}
implies that \eqref{lscond} must be satisfied after at most $\max\{n_k, \,\hat{n}_k\}+1$ inner iterations, which leads to a contradiction.
\endproof

Now, we are ready to prove our main convergence result, which characterizes a cluster point of the sequence generated by Algorithm \ref{alg_NAUM}. Our proof of statement (ii) in the following theorem is similar to that of \cite[Lemma 4]{wnf2009sparse}. However, the arguments involved are more intricate since we have two blocks of variables in our line search loop.

\begin{theorem}\label{subconvergence}
Suppose that Assumption \ref{assumfun} holds. Let $\{(X^k, Y^k)\}$ be the sequence generated by Algorithm \ref{alg_NAUM}. Then,
\begin{itemize}
\item[{\rm (i)}] \textbf{(boundedness of sequence)} the sequences $\{(X^k, Y^k)\}$, $\{\bar{\mu}_k\}$ and $\{\bar{\sigma}_k\}$ are bounded;

\item[{\rm (ii)}] \textbf{(diminishing successive changes)} $\lim_{k\rightarrow \infty} \|X^{k+1} - X^k\|_F + \|Y^{k+1} - Y^k\|_F  = 0$;

\item[{\rm (iii)}] \textbf{(global subsequential convergence)} any cluster point $(X^*, Y^*)$ of $\{(X^k, Y^k)\}$ is a stationary point of $\mathcal{F}$.
\end{itemize}
\end{theorem}
\beginproof
\textit{Statement (i)}. We first show that
\begin{eqnarray}\label{bdbyin}
\mathcal{F} (X^k, Y^k) \leq \mathcal{F} (X^0, Y^0)
\end{eqnarray}
for all $k \geq 1$. We will prove it by induction. Indeed, for $k=1$, it follows from Proposition \ref{guals} that
\begin{eqnarray*}
\mathcal{F}(X^1, Y^1) - \mathcal{F}(X^0, Y^0) \leq-\frac{c}{2}\left(\|X^1-X^0\|_F^2+\|Y^1-Y^0\|_F^2\right)\leq 0
\end{eqnarray*}
is satisfied after finitely many inner iterations. Hence, \eqref{bdbyin} holds for $k=1$. We now suppose that \eqref{bdbyin} holds for all $k \leq K$ for some integer $K \geq 1$. Then, we only need to show that \eqref{bdbyin} also holds for $k=K+1$. For $k=K+1$, we have
\begin{eqnarray*}
\begin{aligned}
&\mathcal{F}(X^{K+1}, Y^{K+1}) - \mathcal{F}(X^0, Y^0)
\leq \mathcal{F}(X^{K+1}, Y^{K+1}) - \max\limits_{[K-N]_{+}\leq i\leq K}\mathcal{F}(X^i, Y^i)   \\
&\qquad \leq-\frac{c}{2}\left(\|X^{K+1}-X^{K}\|_F^2+\|Y^{K+1}-Y^{K}\|_F^2\right) \leq 0,
\end{aligned}
\end{eqnarray*}
where the first inequality follows from the induction hypothesis and the second inequality follows from \eqref{lscond}. Hence, \eqref{bdbyin} holds for $k=K+1$. This completes the induction. Then, from \eqref{bdbyin}, we have that for any $k \geq 0$,
\begin{eqnarray*}
\mathcal{F}(X^0, Y^0) \geq \mathcal{F}(X^k, Y^k)  =\Psi(X^k)+\Phi(Y^k)+\frac{1}{2}\left\|\mathcal{A}(X^k(Y^k)^{\top})-\bm{b}\right\|^2,
\end{eqnarray*}
which, together with (a1) in Assumption \ref{assumfun}, implies that $\{X^k\}$, $\{Y^k\}$ and $\{\|\mathcal{A}(X^k(Y^k)^{\top})-\bm{b}\|\}$ are bounded. Moreover, from Step 2 and Step 3 in Algorithm \ref{alg_NAUM}, it is easy to see $\bar{\mu}_k\leq\mu_k^{\max}=(\alpha+2\gamma\rho)\|Y^k\|^2+c$ for all $k$. Since $\{Y^k\}$ is bounded, the sequences $\{\mu_k^{\max}\}$ and $\{\bar{\mu}_k\}$ are bounded. Next, we prove the boundedness of $\{\bar{\sigma}_k\}$. Indeed, at the $k$-th iteration, there are three possibilities:
\begin{itemize}
\item $\bar{\mu}_k<\mu_k^{\max}$: In this case, we have $\bar{\sigma}_k\leq\sigma^0_k \tau^{\tilde{n}_k}\leq\sigma^{\max} \tau^{\tilde{n}_k}$, where $\tilde{n}_k$ denotes the number of inner iterations for the line search at the $k$-th iteration and $\tilde{n}_k \leq \max\left\{1, \,\left\lfloor\frac{\log(\mu_k^{\max})-\log(\mu^{\min})}{\log\tau} + 2 \right\rfloor\right\}$ (see \eqref{upbnk} and the discussions preceding it).

\item $\bar{\mu}_k=\mu_k^{\max}$ and $\bar{\sigma}_k > \sigma_k^{\max}$: In this case, we have $\bar{\sigma}_k\leq\sigma^0_k \tau^{\tilde{n}_k}\leq\sigma^{\max} \tau^{\tilde{n}_k}$, where $\tilde{n}_k \leq \max\Big{\{}1$, $\left\lfloor\frac{\log(\mu_k^{\max})-\log(\mu^{\min})}{\log\tau} + 2 \right\rfloor\Big{\}}$.

\item Otherwise, we have $\bar{\sigma}_k \leq \sigma_k^{\max}=(\alpha+2\gamma\rho)\|X^{k+1}\|^2+c$.
\end{itemize}
Note that $\{\tilde{n}_k\}$ is bounded as $\{\mu_k^{\max}\}$ is bounded. Thus, $\{\bar{\sigma}_k\}$ is bounded as the sequences $\{X^k\}$ and $\{\tilde{n}_k\}$ are bounded. This proves statement (i).

\textit{Statement (ii)}. We first claim that any cluster point of $\{(X^k, Y^k)\}$ is in $\mathrm{dom}\mathcal{F}$. Since $\{(X^k, Y^k)\}$ is bounded from statement (i), there exists at least one cluster point. Suppose that $(X^*, Y^*)$ is a cluster point of $\{(X^k, Y^k)\}$ and let $\{(X^{k_i}, \!Y^{k_i})\}$ be a convergent subsequence such that $\lim\limits_{i\rightarrow\infty}\!(X^{k_i}, Y^{k_i}) =(X^*, Y^*)$. Then, from the lower semicontinuity of $\mathcal{F}$ (since $\Psi$, $\Phi$ are closed by {\bf(a1)} in Assumption \ref{assumfun}) and \eqref{bdbyin}, we have
\begin{eqnarray*}
\mathcal{F}(X^*, Y^*)\leq\lim\limits_{i\rightarrow\infty}\!\mathcal{F}(X^{k_i}, Y^{k_i})\leq\mathcal{F}(X^0, Y^0),
\end{eqnarray*}
which implies that $\mathcal{F}(X^*, Y^*)$ is finite and hence $(X^*, Y^*) \in \mathrm{dom}\mathcal{F}$.

For notational simplicity, from now on, we let $\Delta_{X^{k}}:=X^{k+1}-X^{k}$, $\Delta_{Y^{k}}:=Y^{k+1}-Y^{k}$, $\Delta_{Z^{k}}:=Z^{k+1}-Z^{k}$ and
\begin{eqnarray}\label{defsk}
\ell(k) = \arg\max\limits_{i}\{\,\mathcal{F}(X^i, Y^i)\,:\,i = [k-N]_+, \cdots, k\,\}.
\end{eqnarray}
Then, the line search criterion \eqref{lscond} can be rewritten as
\begin{eqnarray}\label{lscond2}
\mathcal{F}(X^{k+1}, Y^{k+1}) - \mathcal{F}(X^{\ell(k)}, Y^{\ell(k)})\leq-\frac{c}{2}\left(\|\Delta_{X^{k}}\|_F^2+\|\Delta_{Y^{k}}\|_F^2\right)\leq 0.
\end{eqnarray}
Observe that
\begin{eqnarray*}
\begin{aligned}
&\quad\mathcal{F}(X^{\ell(k+1)}, \,Y^{\ell(k+1)})
=\max\limits_{[k+1-N]_{+}\leq i\leq k+1}\mathcal{F}(X^i, \,Y^i) \\
&=\max\left\{\!\mathcal{F}(X^{k+1}, Y^{k+1}), \max\limits_{[k+1-N]_{+}\leq i\leq {k}}\mathcal{F}(X^i, Y^i)\!\right\}
\stackrel{(\mathrm{i})}{\leq}\max\left\{\!\mathcal{F}(X^{\ell(k)}, Y^{\ell(k)}), \max\limits_{[k+1-N]_{+}\leq i\leq {k}}\mathcal{F}(X^i, Y^i)\!\right\}  \\
&\leq\max\left\{\mathcal{F}(X^{\ell(k)}, Y^{\ell(k)}), \max\limits_{[k-N]_{+}\leq i\leq {k}}\mathcal{F}(X^i, Y^i)\right\} \stackrel{(\mathrm{ii})}{=}\max\left\{\mathcal{F}(X^{\ell(k)}, Y^{\ell(k)}), \,\mathcal{F}(X^{\ell(k)}, Y^{\ell(k)})\right\}  \\
&=\mathcal{F}(X^{\ell(k)}, \,Y^{\ell(k)}),
\end{aligned}
\end{eqnarray*}
where (i) follows from \eqref{lscond2} and (ii) follows from \eqref{defsk}. The sequence $\{\mathcal{F}(X^{\ell(k)}, \,Y^{\ell(k)})\}$ is therefore non-increasing. Since $\mathcal{F}(X^{\ell(k)}, \,Y^{\ell(k)})$ is also bounded from below (due to (a1) in Assumption \ref{assumfun}), we conclude that there exists a number $\widetilde{\mathcal{F}}$ such that
\begin{eqnarray}\label{limsk}
\lim\limits_{k\rightarrow\infty}\mathcal{F}(X^{\ell(k)}, \,Y^{\ell(k)}) = \widetilde{\mathcal{F}}.
\end{eqnarray}

We next prove by induction that for all $j\geq1$,
\begin{numcases}{}
\lim\limits_{k\rightarrow\infty} \Delta_{X^{\ell(k)-j}}=\lim\limits_{k\rightarrow\infty} \Delta_{Y^{\ell(k)-j}} = 0,  \label{keylim1} \\
\lim\limits_{k\rightarrow\infty}\mathcal{F}(X^{\ell(k)-j}, Y^{\ell(k)-j}) = \widetilde{\mathcal{F}} \label{keylim2}.
\end{numcases}
We first prove \eqref{keylim1} and \eqref{keylim2} for $j=1$. Applying \eqref{lscond2} with $k$ replaced by $\ell(k)-1$, we obtain
\begin{eqnarray*}
\mathcal{F}(X^{\ell(k)}, Y^{\ell(k)}) - \mathcal{F}(X^{\ell(\ell(k)-1)}, Y^{\ell(\ell(k)-1)})
\leq-\frac{c}{2}\left(\|\Delta_{X^{\ell(k)-1}}\|_F^2+\|\Delta_{Y^{\ell(k)-1}}\|_F^2\right).
\end{eqnarray*}
Thus, from this and \eqref{limsk}, we have
\begin{eqnarray*}
\lim\limits_{k\rightarrow\infty}~\frac{c}{2}\left(\|\Delta_{X^{\ell(k)-1}}\|_F^2
+\|\Delta_{Y^{\ell(k)-1}}\|_F^2\right) = 0,
\end{eqnarray*}
which implies that
\begin{eqnarray}\label{limsk1}
\lim\limits_{k\rightarrow\infty} \Delta_{X^{\ell(k)-1}}=\lim\limits_{k\rightarrow\infty} \Delta_{Y^{\ell(k)-1}} = 0.
\end{eqnarray}
Then, from \eqref{limsk} and \eqref{limsk1}, we have
\begin{eqnarray*}
\begin{aligned}
\widetilde{\mathcal{F}}~
&=~\lim\limits_{k\rightarrow\infty}\mathcal{F}(X^{\ell(k)}, Y^{\ell(k)}) \\
&=~\lim\limits_{k\rightarrow\infty}\mathcal{F}(X^{\ell(k)-1}+\Delta_{X^{\ell(k)-1}}, Y^{\ell(k)-1}+\Delta_{Y^{\ell(k)-1}}) \\
&=~\lim\limits_{k\rightarrow\infty}\mathcal{F}(X^{\ell(k)-1}, Y^{\ell(k)-1}),
\end{aligned}
\end{eqnarray*}
where the last equality follows because $\{(X^k, Y^k)\}$ is bounded, any cluster point of $\{(X^k, Y^k)\}$ is in $\mathrm{dom}\mathcal{F}$ and $\mathcal{F}$ is uniformly continuous on any compact subset of $\mathrm{dom}\mathcal{F}$ under (a1) in Assumption \ref{assumfun}. Thus, \eqref{keylim1} and \eqref{keylim2} hold for $j=1$.

We next suppose that \eqref{keylim1} and \eqref{keylim2} hold for $j=J$ for some $J \geq 1$. It remains to show that they also hold for $j = J+1$. Indeed, from \eqref{lscond2} with $k$ replaced by $\ell(k)-J-1$ (here, without loss of generality, we assume that $k$ is large enough such that $\ell(k)-J-1$ is nonnegative), we have
\begin{eqnarray*}
\mathcal{F}(X^{\ell(k)-J}, Y^{\ell(k)-J}) - \mathcal{F}(X^{\ell(\ell(k)-J-1)}, Y^{\ell(\ell(k)-J-1)})
\leq-\frac{c}{2}\left(\|\Delta_{X^{\ell(k)-J-1}}\|_F^2+\|\Delta_{Y^{\ell(k)-J-1}}\|_F^2\right),
\end{eqnarray*}
which implies that
\begin{eqnarray*}
\|\Delta_{X^{\ell(k)-J-1}}\|_F^2+\|\Delta_{Y^{\ell(k)-J-1}}\|_F^2
\leq \frac{2}{c}\left(\mathcal{F}(X^{\ell(\ell(k)-J-1)}, Y^{\ell(\ell(k)-J-1)})-\mathcal{F}(X^{\ell(k)-J}, Y^{\ell(k)-J})\right).
\end{eqnarray*}
This together with \eqref{limsk} and the induction hypothesis implies that
\begin{eqnarray*}
\lim\limits_{k\rightarrow\infty} \Delta_{X^{\ell(k)-(J+1)}}=\lim\limits_{k\rightarrow\infty} \Delta_{Y^{\ell(k)-(J+1)}} = 0.
\end{eqnarray*}
Thus, \eqref{keylim1} holds for $j=J+1$. From this, we further have
\begin{eqnarray*}
\begin{aligned}
\lim\limits_{k\rightarrow\infty}\mathcal{F}(X^{\ell(k)-(J+1)}, \,Y^{\ell(k)-(J+1)})
&=\lim\limits_{k\rightarrow\infty}\mathcal{F}(X^{\ell(k)-J}-\Delta_{X^{\ell(k)-(J+1)}}, \,Y^{\ell(k)-J}-\Delta_{Y^{\ell(k)-(J+1)}}) \\
&=\lim\limits_{k\rightarrow\infty}\mathcal{F}(X^{\ell(k)-J}, \,Y^{\ell(k)-J}) =\widetilde{\mathcal{F}},
\end{aligned}
\end{eqnarray*}
where the second equality follows because $\{(X^k, Y^k)\}$ is bounded, any cluster point of $\{(X^k, Y^k)\}$ is in $\mathrm{dom}\mathcal{F}$ and $\mathcal{F}$ is uniformly continuous on any compact subset of $\mathrm{dom}\mathcal{F}$ under (a1) in Assumption \ref{assumfun}. Hence, \eqref{keylim2} also holds for $j=J+1$. This completes the induction.

We are now ready to prove the main result in this statement. Indeed, from \eqref{defsk}, we can see $k-N\leq \ell(k)\leq k$ (without loss of generality, we assume that $k$ is large enough such that $k \geq N$). Thus, for any $k$, we must have $k-N-1=\ell(k)-j_k$ for $1\leq j_k \leq N+1$. Then, we have
\begin{eqnarray*}
\begin{aligned}
&\|\Delta_{X^{k-N-1}}\|_F=\|\Delta_{X^{\ell(k)-j_k}}\|_F  \leq \max\limits_{1\leq j \leq N+1} \|\Delta_{X^{\ell(k)-j}}\|_F, \\
&\|\Delta_{Y^{k-N-1}}\|_F=\|\Delta_{Y^{\ell(k)-j_k}}\|_F  \leq \max\limits_{1\leq j \leq N+1} \|\Delta_{Y^{\ell(k)-j}}\|_F.
\end{aligned}
\end{eqnarray*}
This together with \eqref{keylim1} implies that
\begin{eqnarray*}
\begin{aligned}
&\lim\limits_{k\rightarrow\infty} \Delta_{X^{k}} = \lim\limits_{k\rightarrow\infty} \Delta_{X^{k-N-1}} = 0,  \\
&\lim\limits_{k\rightarrow\infty} \Delta_{Y^{k}} = \lim\limits_{k\rightarrow\infty} \Delta_{Y^{k-N-1}} = 0.  \\
\end{aligned}
\end{eqnarray*}
This proves the statement (ii).

\textit{Statement (iii)}. Again, let $(X^*, Y^*)$ be a cluster point of $\{(X^k, Y^k)\}$ and let $\{(X^{k_i}, \!Y^{k_i})\}$ be a convergent subsequence such that $\lim\limits_{i\rightarrow\infty}\!(X^{k_i}, Y^{k_i}) =(X^*, Y^*)$. Recall that $(X^*, Y^*) \in \mathrm{dom}\mathcal{F}$. On the other hand, it is easy to see from \eqref{Zupdate} that $\lim\limits_{i\rightarrow\infty}\,Z^{k_i} = Z^*$, where $Z^*$ is given by \eqref{addZ}. Thus, it can be shown as in \eqref{zequa} that
\begin{eqnarray}\label{optz}
\alpha(Z^* - X^*(Y^*)^{\top}) + \beta\mathcal{A}^*(\mathcal{A}(Z^*)-\bm{b}) = 0.
\end{eqnarray}
We next show that
\begin{numcases}{}
0 \in \partial\Psi(X^*) + \alpha(X^*(Y^*)^{\top}-Z^*)Y^*,   \label{optx}\\
0 \in \partial\Phi(Y^*) + \alpha(X^*(Y^*)^{\top}-Z^*)^{\top}X^*. \label{opty}
\end{numcases}
We start by showing \eqref{optx} in the following cases:
\begin{itemize}
\item Proximal~\&~Prox-linear: In these two cases, passing to the limit along $\{(X^{k_i}, Y^{k_i})\}$ in \eqref{Xoptcon} or \eqref{Xlinoptcon} with $X^{k_i+1}$ in place of $U$ and $\bar{\mu}_{k_i}$ in place of $\mu_{k}$, and invoking (a1) in Assumption \ref{assumfun}, statements (i), (ii), $(X^*, Y^*) \in \mathrm{dom}\mathcal{F}$ and \eqref{robust}, we obtain \eqref{optx}.

\item Hierarchical-prox: In this case, passing to the limit along $\{(X^{k_i}, Y^{k_i})\}$ in \eqref{Xhioptcon} with $X^{k_i+1}$ in place of $U$ and $\bar{\mu}_{k_i}$ in place of $\mu_{k}$, and invoking (a1) in Assumption \ref{assumfun}, statements (i), (ii), $(X^*, Y^*) \in \mathrm{dom}\mathcal{F}$ and \eqref{robust}, we have
    \begin{eqnarray*}
    0 \in \partial\psi_i(\bm{x}^*_{i}) + \alpha(X^*(Y^*)^{\top}-Z^*)\bm{y}^*_i
    \end{eqnarray*}
    for any $i=1,2,\cdots,r$. Then, stacking them up, we obtain \eqref{optx}.

\end{itemize}
Similarly, we can obtain \eqref{opty}. Thus, combining \eqref{optz}, \eqref{optx} and \eqref{opty}, we see that $(X^*, Y^*, Z^*)$ is a stationary point of $\Theta_{\alpha,\beta}$, which further implies $(X^*, Y^*)$ is a stationary point of $\mathcal{F}$ from Theorem \ref{modelstaeq}. This proves statement (iii).
\endproof

\begin{remark}[\textbf{Comment on (a3) in Assumption \ref{assumfun}}]
If $\Phi$ and $\Psi$ are the indicator functions of some nonempty closed sets, the results in Theorem \ref{subconvergence} remain valid under the weaker condition on $\alpha$ and $\beta$ that $\frac{1}{\alpha}+\frac{1}{\beta}>0$ with a slight modification in \eqref{lscond} of Algorithm \ref{alg_NAUM}. Indeed, when $\Phi$ and $\Psi$ are the indicator functions, one can see from Remark \ref{remark1} and the proofs of Lemma \ref{suffdes} and Proposition \ref{guals} that if $\frac{1}{\alpha}+\frac{1}{\beta}>0$, then
\begin{eqnarray*}
\begin{aligned}
&~~\mathcal{F}(U, V) - \mathcal{F}(X^k, Y^k)
= {\textstyle\left(\frac{1}{\alpha}+\frac{1}{\beta}\right)\left(\Theta_{\alpha,\beta}(U, V, W) - \Theta_{\alpha,\beta}(X^k, Y^k, Z^k)\right)} \\
&\leq-{\textstyle\left(\frac{1}{\alpha}+\frac{1}{\beta}\right)
\left(\frac{\mu_k-(\alpha+2\gamma\rho)\|Y^{k}\|^2}{2}\cdot\|U-X^{k}\|_F^2
+\frac{\sigma_k-(\alpha+2\gamma\rho)\|U\|^2}{2}\cdot\|V-Y^{k}\|_F^2\right)},
\end{aligned}
\end{eqnarray*}
and the line search criterion is well defined with $c$ replaced by $\left(\frac{1}{\alpha}+\frac{1}{\beta}\right)c$. Moreover, recalling \cite[Exercise 8.14]{rw1998variational}, we see that $\partial\Psi$ and $\partial\Phi$ are normal cones. Thus, following Remark \ref{conecond} and the similar augments in Theorem \ref{subconvergence}, we can obtain the same results when $\frac{1}{\alpha}+\frac{1}{\beta}>0$ with $c$ replaced by $\left(\frac{1}{\alpha}+\frac{1}{\beta}\right)c$ in \eqref{lscond} of Algorithm \ref{alg_NAUM}.
\end{remark}

\begin{remark}[\textbf{Comments on updating $\mu_k^{\max}$ and $\sigma_k^{\max}$}]
In Algorithm \ref{alg_NAUM}, we need to evaluate
$\mu_k^{\max}=(\alpha+2\gamma\rho)\|Y^k\|^2+c$ and $\sigma_k^{\max}=(\alpha+2\gamma\rho)\|U\|^2+c$ in each iteration. However, computing the spectral norms of $Y^k$ and $U$ might be costly, especially when $r$ is large. Hence, in our experiments, instead of computing $\|Y^k\|^2$ and $\|U\|^2$, we compute $\|Y^k\|^2_F$ and $\|U\|^2_F$, and update $\mu_k^{\max}$ and $\sigma_k^{\max}$ by $\mu_k^{\max}=(\alpha+2\gamma\rho)\|Y^k\|_F^2+c$ and $\sigma_k^{\max}=(\alpha+2\gamma\rho)\|U\|_F^2+c$ instead. Since $\|Y^k\|\leq\|Y^k\|_F$ and $\|U\|\leq\|U\|_F$, it follows from \eqref{succhange} that
\begin{eqnarray*}
\mathcal{F}(U, V) - \mathcal{F}(X^k, Y^k)
\leq-{\textstyle\frac{\mu_k-(\alpha+2\gamma\rho)\|Y^{k}\|_F^2}{2}}\,\|U-X^{k}\|_F^2
-{\textstyle\frac{\sigma_k-(\alpha+2\gamma\rho)\|U\|_F^2}{2}}\,\|V-Y^{k}\|_F^2.
\end{eqnarray*}
Then, one can show that Proposition \ref{guals} and Theorem \ref{subconvergence} remain valid. In addition, we compute the quantities $\|U\|^2_F$ and $\|Y^k\|^2_F$ by $\mathrm{tr}(U^{\top}U)$ and $\mathrm{tr}((Y^k)^{\top}Y^k)$, respectively. For some cases, the matrices $U^{\top}U$ and $(Y^k)^{\top}Y^k$ can be used repeatedly in updating the variables and evaluating the objective value and successive changes to reduce the cost of line search; see a concrete example in Section \ref{secNMF}.
\end{remark}

\section{Numerical experiments}\label{secnum}

In this section, we conduct numerical experiments to test our algorithm for NMF and MC on real datasets. All experiments are run in MATLAB R2015b on a 64-bit PC with an Intel Core i7-4790 CPU (3.60 GHz) and 32 GB of RAM equipped with Windows 10 OS.

\subsection{Non-negative matrix factorization}\label{secNMF}

We first consider NMF
\begin{eqnarray}\label{nmfmodelbasic}
\min \limits_{X,Y}~~\frac{1}{2}\left\|XY^{\top}-M\right\|_F^2 \quad \mathrm{s.t.} \quad X \geq 0, ~~Y \geq 0,
\end{eqnarray}
where $X\in\mathbb{R}^{m\times r}$ and $Y\in\mathbb{R}^{n\times r}$ are decision variables. Note that the feasible set of \eqref{nmfmodelbasic} is unbounded. We hence focus on the following model:
\begin{eqnarray}\label{nmfmodel}
\min \limits_{X,Y}~~\frac{1}{2}\left\|XY^{\top}-M\right\|_F^2 \quad \mathrm{s.t.} \quad 0 \leq X \leq X^{\max}, ~0 \leq Y \leq Y^{\max},
\end{eqnarray}
where $X^{\max} \geq 0$ and $Y^{\max} \geq 0$ are upper bound matrices. One can show that, when $X^{\max}_{ij}$ and $Y^{\max}_{ij}$ are sufficiently large\footnote{The estimations of $X^{\max}_{ij}$ and $Y^{\max}_{ij}$ have been discussed in \cite[Page 67]{g2011nonnegative}.}, solving \eqref{nmfmodel} gives a solution of \eqref{nmfmodelbasic}. In our experiments, for simplicity, we set $X^{\max}_{ij}=10^{16}$ and $Y^{\max}_{ij}=10^{16}$ for all $(i, j)$. Now, we see that \eqref{nmfmodel} corresponds to \eqref{MFPmodel} with $\Psi(X)=\delta_{\mathcal{X}}(X)$, $\Phi(Y)=\delta_{\mathcal{Y}}(Y)$ and $\mathcal{A}=\mathcal{I}$, where $\mathcal{X}=\{X\in\mathbb{R}^{m\times r}: 0 \leq X \leq X^{\max}\}$ and $\mathcal{Y}=\{Y\in\mathbb{R}^{n\times r}: 0 \leq Y \leq Y^{\max}\}$. We apply NAUM to solving \eqref{nmfmodel}, and use \eqref{Xhiupdate} and \eqref{Yhiupdate} to update $U$ and $V$. The specific updates of $Z^k$, $\bm{u}_i$ and $\bm{v}_i$ are
\begin{eqnarray*}
\begin{aligned}
Z^{k} &= {\textstyle\frac{\alpha}{\alpha+\beta}}X^k(Y^k)^{\top} + {\textstyle\frac{\beta}{\alpha+\beta}}M,  \\
\bm{u}_i &= \max\left\{0, ~\min\left\{\bm{x}^{\max}_i,~\frac{\alpha P_i^k\bm{y}^k_i+\mu_k\bm{x}^k_i}{\alpha\|\bm{y}_i^k\|^2+\mu_k}\right\}\right\},~~i=1,2\cdots,r, \\
\bm{v}_i &= \max\left\{0,~\min\left\{\bm{y}^{\max}_i, ~\frac{\alpha (Q_i^k)^{\top}\bm{u}_i+\sigma_k\bm{y}^k_i}{\alpha\|\bm{u}_i\|^2+\sigma_k}\right\}\right\},~~i=1,2\cdots,r, \\
\end{aligned}
\end{eqnarray*}
where $P^k_i$ and $Q^k_i$ are defined in \eqref{defPQ}. Note that here it is not necessary to update $Z^k$ explicitly. Indeed, we can directly compute $P_i^k\bm{y}^k_i$ and $(Q_i^k)^{\top}\bm{u}_i$ by substituting $Z^k$ as below:
\begin{eqnarray}\label{substz1}
\begin{aligned}
&P_i^k\bm{y}^k_i \!=\! {\textstyle\frac{\alpha}{\alpha+\beta}}X^k(Y^k)^{\top}\!\bm{y}^k_i
\!+\!{\textstyle\frac{\beta}{\alpha+\beta}}M\bm{y}^k_i
\!-\!{\textstyle\sum_{j=1}^{i-1}}\bm{u}_j(\bm{y}^k_j)^{\top}\!\bm{y}^k_i
\!-\!{\textstyle\sum_{j=i+1}^r}\bm{x}^{k}_j(\bm{y}^k_j)^{\top}\!\bm{y}^k_i, \\
&(Q_i^k)^{\top}\!\bm{u}_i \!=\! {\textstyle\frac{\alpha}{\alpha+\beta}}Y^k(X^k)^{\top}\!\bm{u}_i
\!+\!{\textstyle\frac{\beta}{\alpha+\beta}}M^{\top}\!\bm{u}_i
\!-\!{\textstyle\sum^{i-1}_{j=1}}\bm{v}_j\bm{u}_j^{\top}\!\bm{u}_i
\!-\!{\textstyle\sum^r_{j=i+1}}\bm{y}^k_j\bm{u}_j^{\top}\!\bm{u}_i.
\end{aligned}
\end{eqnarray}
When computing $X^k(Y^k)^{\top}\!\bm{y}^k_i$ and $Y^k(X^k)^{\top}\!\bm{u}_i$ in the above, we first compute $(Y^k)^{\top}\!\bm{y}^k_i$ and $(X^k)^{\top}\bm{u}_i$ to avoid forming the huge ($m\times n$) matrix $X^k(Y^k)^{\top}$. Moreover, the matrices $(X^k)^{\top}U$, $U^{\top}U$, $(Y^{k})^{\top}Y^k$ and $M^{\top}U$ that have been computed in \eqref{substz1} can be used again to evaluate the successive changes and the objective value as follows:
\begin{eqnarray*}
\begin{aligned}
\|U-X^k\|_F^2 &= \mathrm{tr}(U^{\top}U) - 2\mathrm{tr}((X^k)^{\top}U) + \mathrm{tr}((X^k)^{\top}X^k), \\
\|V-Y^k\|_F^2 &= \mathrm{tr}(V^{\top}V) - 2\mathrm{tr}((Y^k)^{\top}V) + \mathrm{tr}((Y^k)^{\top}Y^k), \\
\|UV^{\top}-M\|_F^2 &= \mathrm{tr}((U^{\top}U)(V^{\top}V)) - 2\mathrm{tr}((M^{\top}U)V^{\top}) + \|M\|_F^2.
\end{aligned}
\end{eqnarray*}
In the above relations, $(X^k)^{\top}X^k$ and $(Y^k)^{\top}Y^k$ can be obtained from $U^{\top}U$ and $V^{\top}V$ in the previous iteration, respectively, and $\|M\|_F^2$ can be computed in advance. Additionally, as we discussed in Remark \ref{upbdpara}, $\mathrm{tr}((Y^k)^{\top}Y^k)$ and $\mathrm{tr}(U^{\top}U)$ can also be used in computing $\mu_k^{\max}$ and $\sigma_k^{\max}$, respectively. These techniques were also used in many popular algorithms for NMF to reduce the computational cost (see, for example, \cite{bblpp2007algorithms,g2011nonnegative,g2014the,ls2001algorithms,wz2013nonnegative}).

The experiments are conducted on the face datasets (dense matrices) and the text datasets (sparse matrices). For face datasets, we use CBCL\footnote{Available in \url{http://cbcl.mit.edu/cbcl/software-datasets/FaceData2.html}.}, ORL\footnote{Available in \url{http://www.cl.cam.ac.uk/research/dtg/attarchive/facedatabase.html}.} \cite{sh1994paramterisation} and the extended Yale Face Database B (e-YaleB)\footnote{Available in \url{http://vision.ucsd.edu/~iskwak/ExtYaleDatabase/ExtYaleB.html}.} \cite{lhk2005acquiring} for our test. CBCL contains 2429 images of faces with $19\times19$ pixels, ORL contains 400 images of faces with $112\times92$ pixels, and e-YaleB contains 2414 images of faces with $168\times192$ pixels. In our experiments, for each face dataset, each image is vectorized and stacked as a column of a data matrix $M$ of size $m \times n$. For text datasets, we use three datasets from the CLUTO toolkit\footnote{Available in \url{http://glaros.dtc.umn.edu/gkhome/cluto/cluto/download}.}. The specific values of $m$ and $n$ for each dataset and the values of $r$ used for our tests are summarized in Table \ref{data}.

\begin{table}[ht]
\setlength{\belowcaptionskip}{6pt}
\caption{Real data sets}\label{data}
\centering \tabcolsep 4pt
\small{
\begin{tabular}{|lcccc|lcccc|}
\hline
\multicolumn{5}{|c|}{Face Datasets (dense matrices)} & \multicolumn{5}{c|}{Text Datasets (sparse matrices)}\\
Data  &  Pixels  & $m$ & $n$ & $r$  & Data  &  Sparsity  & $m$ & $n$ & $r$  \\
\hline
CBCL     &  $19\times19$    & 361     & 2429   & 30, 60  &
classic  &  99.92\%         & 7094    & 41681  & 10, 20   \\
ORL      &  $112\times92$   & 10304   & 400    & 30, 60  &
sports   &  99.14\%         & 8580    & 14870  & 10, 20   \\
e-YaleB  &  $168\times192$  & 32256   & 2414   & 30, 60  &
ohscal   &  99.47\%         & 11162   & 11465  & 10, 20    \\
\hline
\end{tabular}}
\end{table}

The parameters in NAUM are set as follows: $\mu^{\min}=\bar{\mu}_{-1}=1$, $\sigma^{\min}=\bar{\sigma}_{-1}=1$, $\sigma^{\max}=10^6$, $\tau=4$, $c=10^{-4}$, $N=3$, $\mu^0_k=\max\left\{0.1\bar{\mu}_{k-1},\,\mu^{\min}\right\}$ and $\sigma^0_k=\min\left\{\max\left\{0.1\bar{\sigma}_{k-1},\,\sigma^{\min}\right\}, \,\sigma^{\max}\right\}$ for any $k\geq0$. Moreover, we set $\beta=\frac{\alpha}{\alpha-1}$, $\gamma=\max\{0, \,-\alpha, \,-(\alpha+\beta)\}$ and $\rho=\max\left\{1,\,\alpha^2/(\alpha+\beta)^2\right\}$ for some given $\alpha$.

We then compare the performances of NAUM with different $\alpha$. In our comparisons, we initialize NAUM with different $\alpha$ at the same random initialization $(X^0, \,Y^0)$\footnote{We use the Matlab commands: \texttt{X0 = max(0, randn(m, r)); Y0 = max(0, randn(n, r)); X0 = X0/norm(X0,'fro')*sqrt(norm(M, 'fro')); Y0 = Y0/norm(Y0,'fro')*sqrt(norm(M, 'fro'));}} and terminate them if one of the following stopping criteria is satisfied:
\begin{itemize}[leftmargin=6mm]
\item $\frac{|\mathcal{F}_{\mathrm{nmf}}^{k} - \mathcal{F}_{\mathrm{nmf}}^{k-1}|}{\mathcal{F}_{\mathrm{nmf}}^{k}+1} \leq 10^{-4}$ holds for 3 consecutive iterations;
\item $\frac{\|X^k-X^{k-1}\|_F+\|Y^k-Y^{k-1}\|_F}{\|X^k\|_F+\|Y^k\|_F+1} \leq 10^{-4}$ holds,
\end{itemize}
where $\mathcal{F}_{\mathrm{nmf}}^k:=\frac{1}{2}\left\|X^k(Y^k)^{\top}-M\right\|_F^2$ denotes the objective value at $(X^k, Y^k)$. Table \ref{comNAUM} presents the results of NAUM with different $\alpha$ for two face datasets (CBCL and ORL) and $r=30,60$. In the table, ``iter" denotes the number of iterations; ``relerr" denotes the relative error $\frac{\|X^*(Y^*)^{\top}-M\|_F}{\|M\|_F}$, where $(X^*, \,Y^*)$ is a terminating point obtained by each NUAM in a trial; ``time" denotes the computational time (in seconds). All the results presented are the average of 10 independent trials. From Table \ref{comNAUM}, we can see that NAUM with a relatively small $\alpha$ (e.g., 0.6 and 0.8) has better numerical performance. However, $\alpha$ cannot be too small. Observe that NAUM with $\alpha=0.5, 0.4, 0.2$ are not competitive and, surprisingly, $\alpha=0.5$ leads to the worst performance. In view of this, we do not choose $\alpha < 0.6$ in our following experiments for NMF.

\begin{table}[ht]
\setlength{\belowcaptionskip}{6pt}
\caption{Comparisons of NAUM with different $\alpha$}\label{comNAUM}
\centering \tabcolsep 6pt
{\small
\begin{tabular}{|crcr|crcr|}
\hline
$\alpha$ & iter & relerr & time & $\alpha$ & iter & relerr & time \\
\hline
\multicolumn{4}{|c|}{CBCL, $r=30$} & \multicolumn{4}{c|}{CBCL, $r=60$}  \\
\hline
 2.0  &  488 & 1.0519e-01 & 1.72 & 2.0  &  626 & 7.4388e-02 & 4.94 \\
 1.1  &  381 & 1.0448e-01 & 1.35 & 1.1  &  555 & 7.3477e-02 & 4.38 \\
 0.8  &  315 & 1.0426e-01 & 1.09 & 0.8  &  511 & 7.2986e-02 & 4.09 \\
 0.6  &  268 & 1.0406e-01 & 0.94 & 0.6  &  419 & 7.2998e-02 & 3.32 \\
 0.5  &  833 & 1.0593e-01 & 4.74 & 0.5  & 1372 & 7.5864e-02 & 19.49 \\
 0.4  &  440 & 1.0489e-01 & 3.05 & 0.4  &  599 & 7.4568e-02 & 10.02 \\
 0.2  &  556 & 1.0674e-01 & 4.18 & 0.2  &  782 & 7.7654e-02 & 14.30 \\
\hline
\multicolumn{4}{|c|}{ORL, $r=30$} & \multicolumn{4}{c|}{ORL, $r=60$}  \\
\hline
 2.0  &  232 & 1.6673e-01 & 3.45 & 2.0  &  277 & 1.4078e-01 & 7.92 \\
 1.1  &  188 & 1.6619e-01 & 2.78 & 1.1  &  210 & 1.4042e-01 & 6.04 \\
 0.8  &  158 & 1.6603e-01 & 2.33 & 0.8  &  182 & 1.4017e-01 & 5.20 \\
 0.6  &  132 & 1.6578e-01 & 2.01 & 0.6  &  156 & 1.3996e-01 & 4.44 \\
 0.5  &  652 & 1.7216e-01 & 15.79 & 0.5  &  695 & 1.4583e-01 & 32.91 \\
 0.4  &  280 & 1.6615e-01 & 7.55 & 0.4  &  353 & 1.4061e-01 & 19.17 \\
 0.2  &  307 & 1.6753e-01 & 8.71 & 0.2  &  358 & 1.4272e-01 & 20.77 \\
 \hline
\end{tabular}
}
\end{table}

We next compare NAUM with two existing efficient algorithms\footnote{Most existing algorithms are directly developed for \eqref{nmfmodelbasic}. However, they need the assumption that the sequence generated is bounded in their convergence analysis. Although this assumption is uncheckable and may fail, these algorithms always work well in practice. Thus, we directly use these algorithms in our comparisons, rather than modifying them for \eqref{nmfmodel}.} for NMF: the hierarchical alternating least squares (HALS) method\footnote{HALS for \eqref{nmfmodelbasic} is given by
{\scriptsize\begin{eqnarray*}
\begin{aligned}
&\bm{x}_i^{k+1}=\max\left\{0, ~\frac{M\bm{y}^k_i-{\textstyle\sum^{i-1}_{j=1}}\bm{x}^{k+1}_j(\bm{y}^k_j)^{\top}\bm{y}^k_i
-{\textstyle\sum^r_{j=i+1}}\bm{x}^{k}_j(\bm{y}^k_j)^{\top}\bm{y}^k_i}{\|\bm{y}_i^k\|^2}\right\},
~~i=1,\cdots,r,   \\
&\bm{y}_i^{k+1}=\max\left\{0, ~\frac{M^{\top}\bm{x}_i^{k+1}-{\textstyle\sum^{i-1}_{j=1}}\bm{y}^{k+1}_j(\bm{x}^{k+1}_j)^{\top}\bm{x}_i^{k+1}
-{\textstyle\sum^r_{j=i+1}}\bm{y}^k_j(\bm{x}^{k+1}_j)^{\top}\bm{x}_i^{k+1}}{\|\bm{x}_i^{k+1}\|^2}\right\},
~~i=1,\cdots,r.
\end{aligned}
\end{eqnarray*}}}
(see, for example, \cite{cza2007hierarchical,g2011nonnegative,g2014the,gg2012accelerated,lz2009fastnmf,llwy2012sparse}) and the block coordinate descent method for NMF (BCD-NMF\footnote{Available at \url{http://www.math.ucla.edu/~wotaoyin/papers/bcu/nmf/index.html}.}) (see Algorithm 2 in Section 3.2 in \cite{xy2013a}).

To better evaluate the performances of different algorithms, we follow \cite{gg2012accelerated} to use an evolution of the objective function value. To define this evolution, we first define
\begin{eqnarray*}
e(k):=\frac{\mathcal{F}^k - \mathcal{F}_{\min}}{\mathcal{F}^0 - \mathcal{F}_{\min}},
\end{eqnarray*}
where $\mathcal{F}^k$ denotes the objective function value obtained by an algorithm at $(X^k, Y^k)$ and $\mathcal{F}_{\min}$ denotes the minimum of the objective function values obtained among \textit{all} algorithms across \textit{all} initializations. We also use $\mathcal{T}(k)$ to denote the total computational time after completing the $k$-th iteration of an algorithm. Thus, $\mathcal{T}(0)=0$ and $\mathcal{T}(k)$ is non-decreasing with respect to $k$. Then, the evolution of the function value obtained from a particular algorithm with respect to time $t$ is defined as
\begin{eqnarray*}
E(t):=\min\left\{e(k):k\in\{i: \mathcal{T}(i) \leq t\}\right\}.
\end{eqnarray*}
One can see that $0 \leq E(t)\leq1$ (since $0 \leq e(k)\leq1$ for all $k$) and $E(t)$ is non-increasing with respect to $t$. $E(t)$ can be considered as a normalized measure of the reduction of the function value with respect to time. For a given matrix $M$ and a positive integer $r$, one can take the average of $E(t)$ over several independent trials with different initializations, and plot the average $E(t)$ within time $t$ for a given algorithm.

\begin{figure}[ht]
\centering
\subfigure[$\mathrm{T}^{\max}=10$]{\includegraphics[width=7cm]{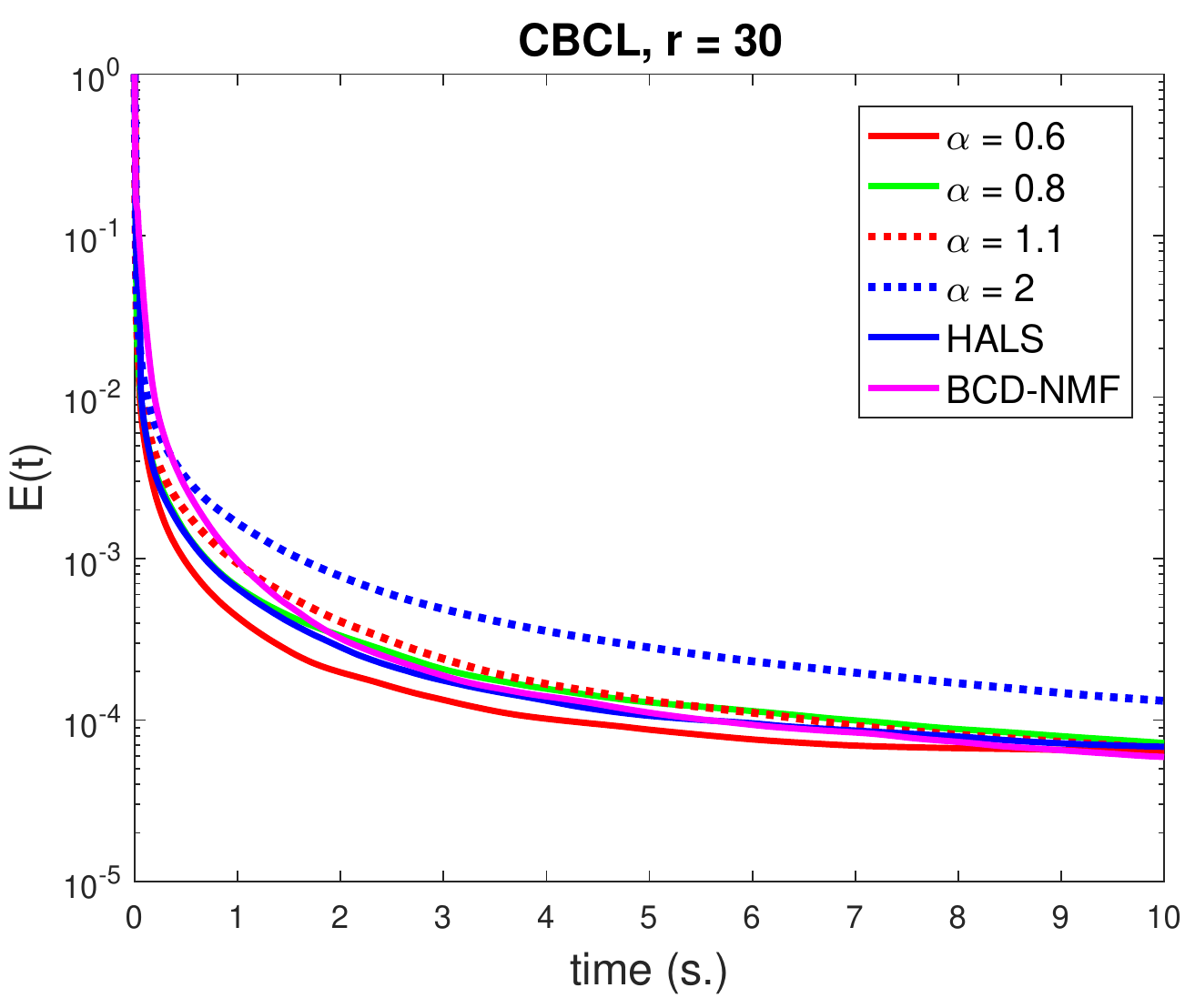}}
\subfigure[$\mathrm{T}^{\max}=20$]{\includegraphics[width=7cm]{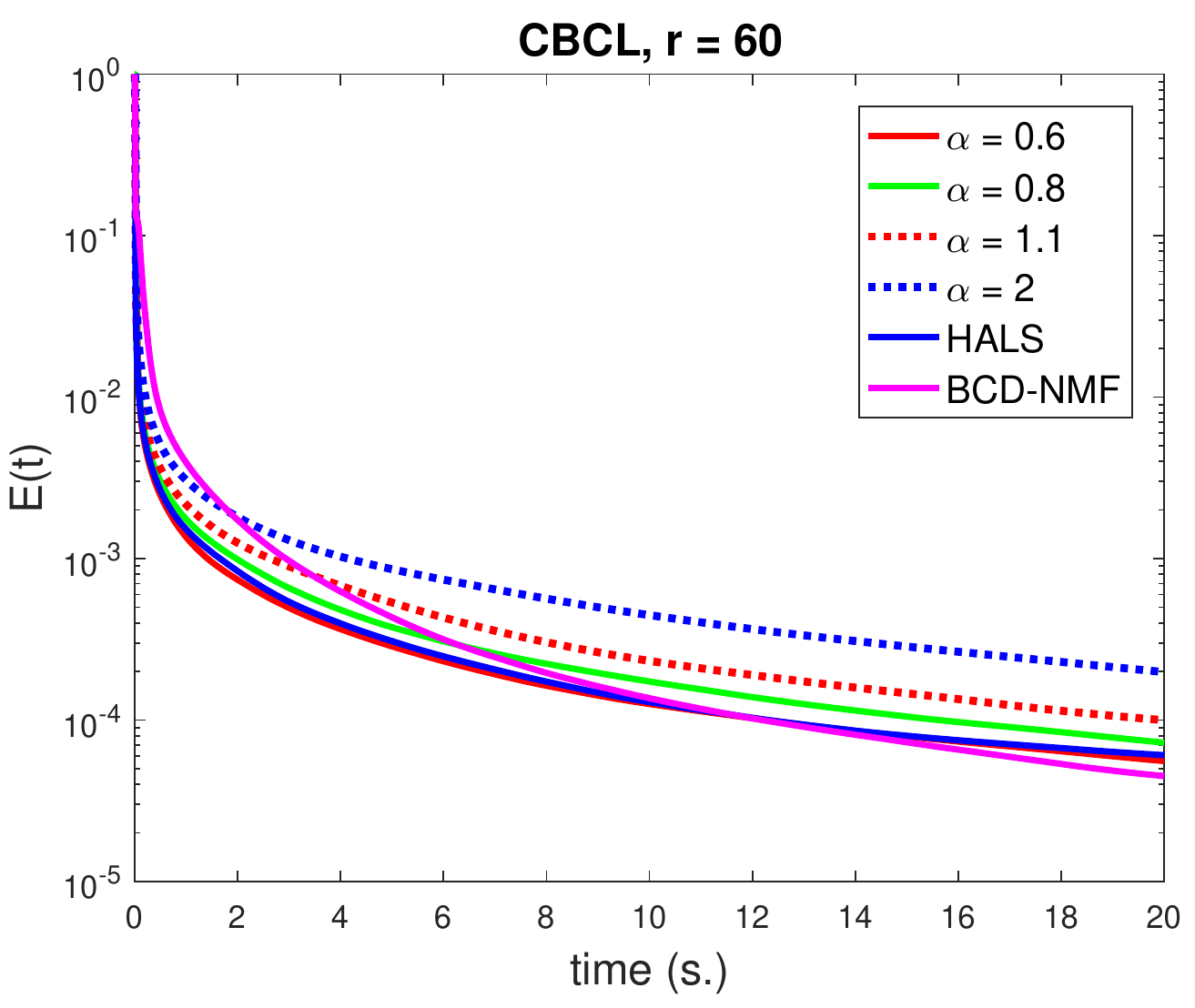}} \\
\subfigure[$\mathrm{T}^{\max}=10$]{\includegraphics[width=7cm]{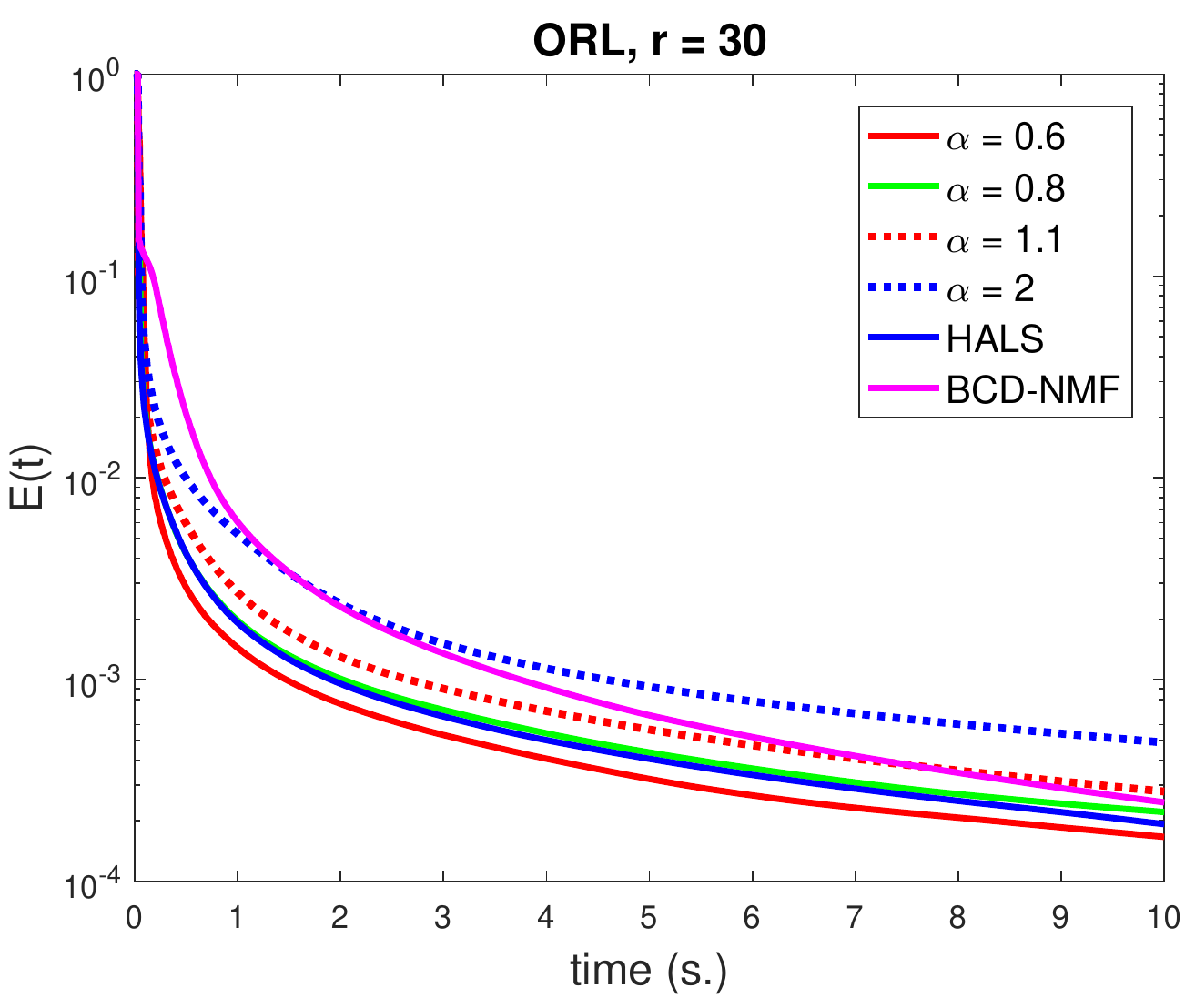}}
\subfigure[$\mathrm{T}^{\max}=20$]{\includegraphics[width=7cm]{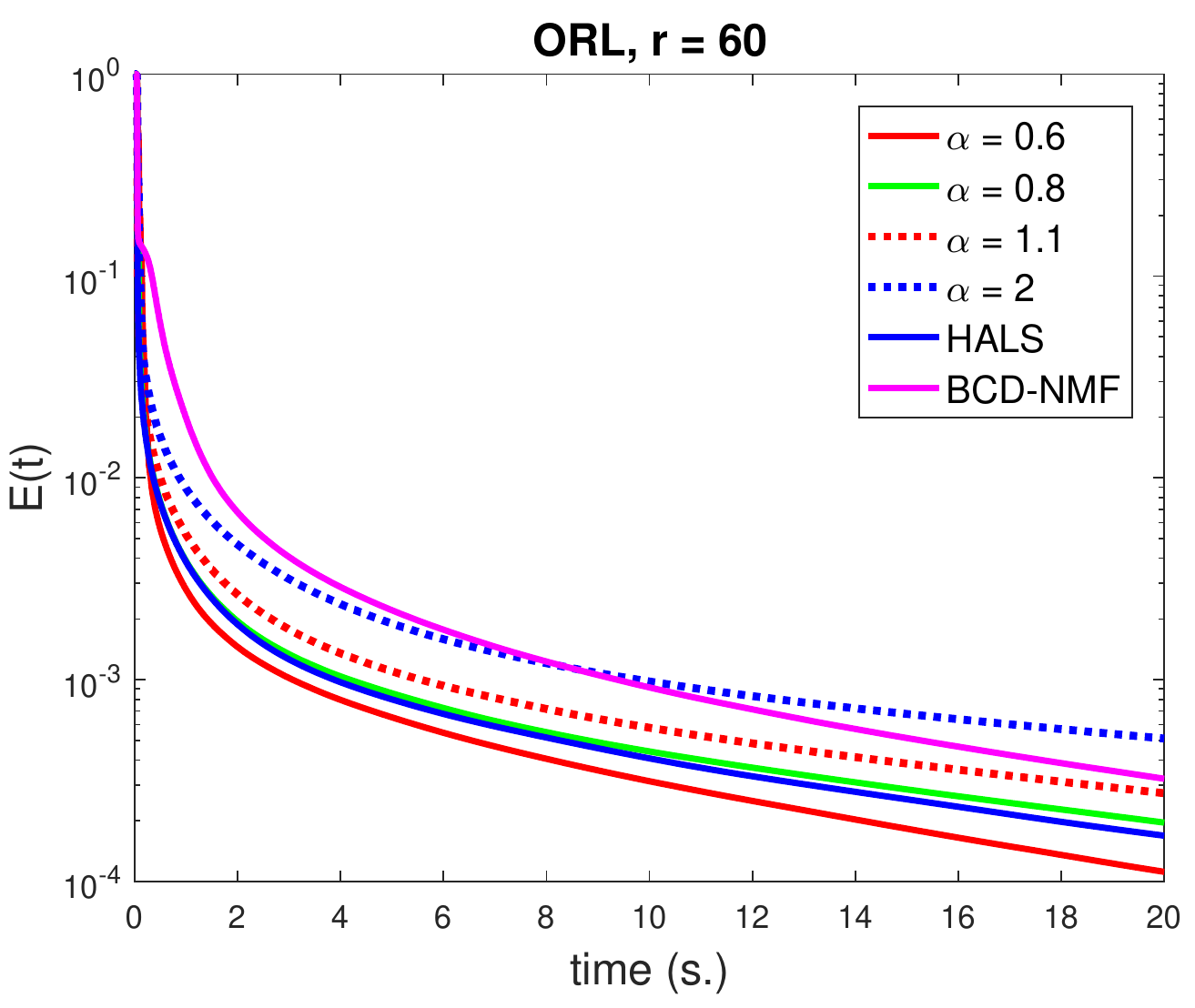}}  \\
\subfigure[$\mathrm{T}^{\max}=60$]{\includegraphics[width=7cm]{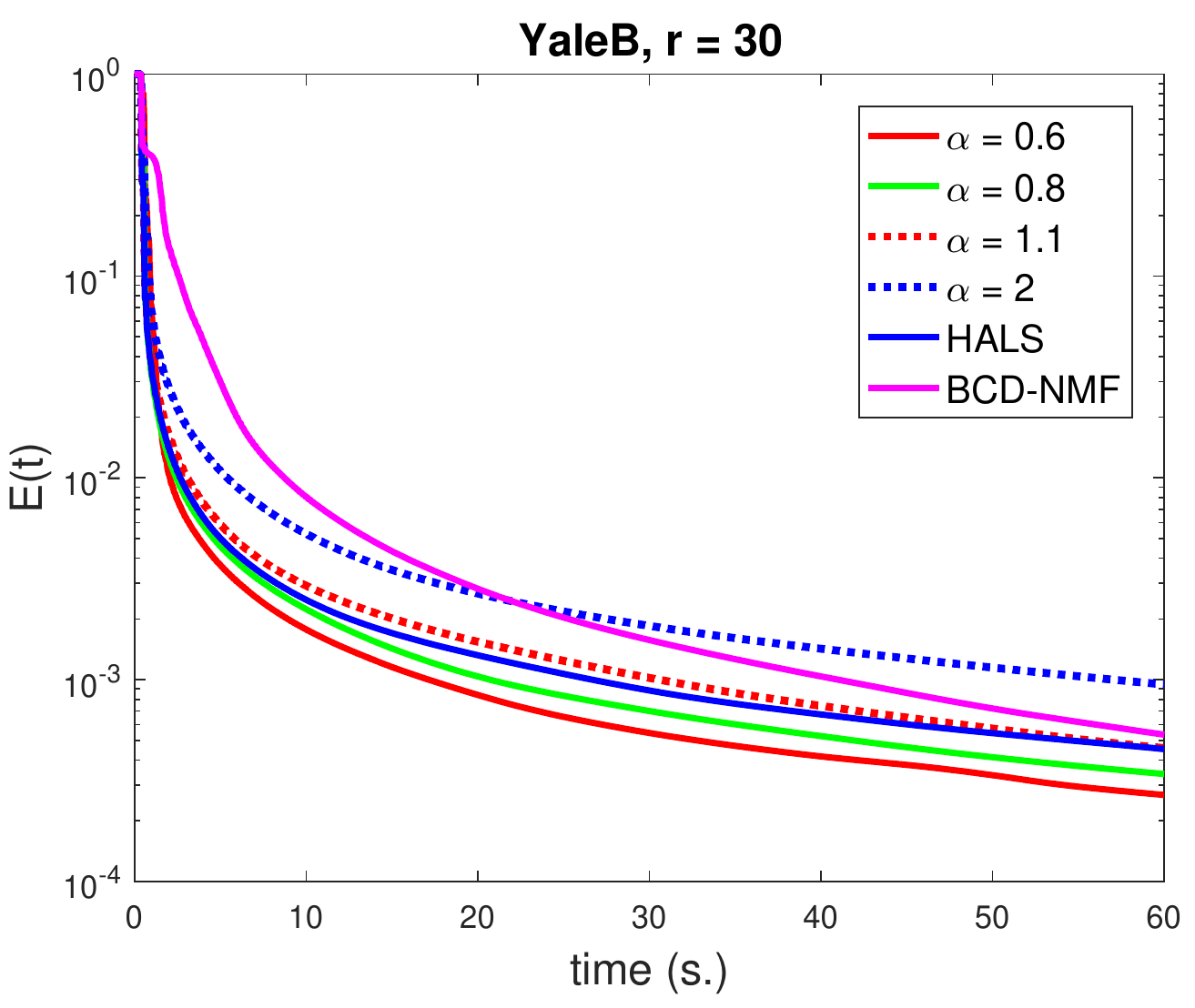}}
\subfigure[$\mathrm{T}^{\max}=120$]{\includegraphics[width=7cm]{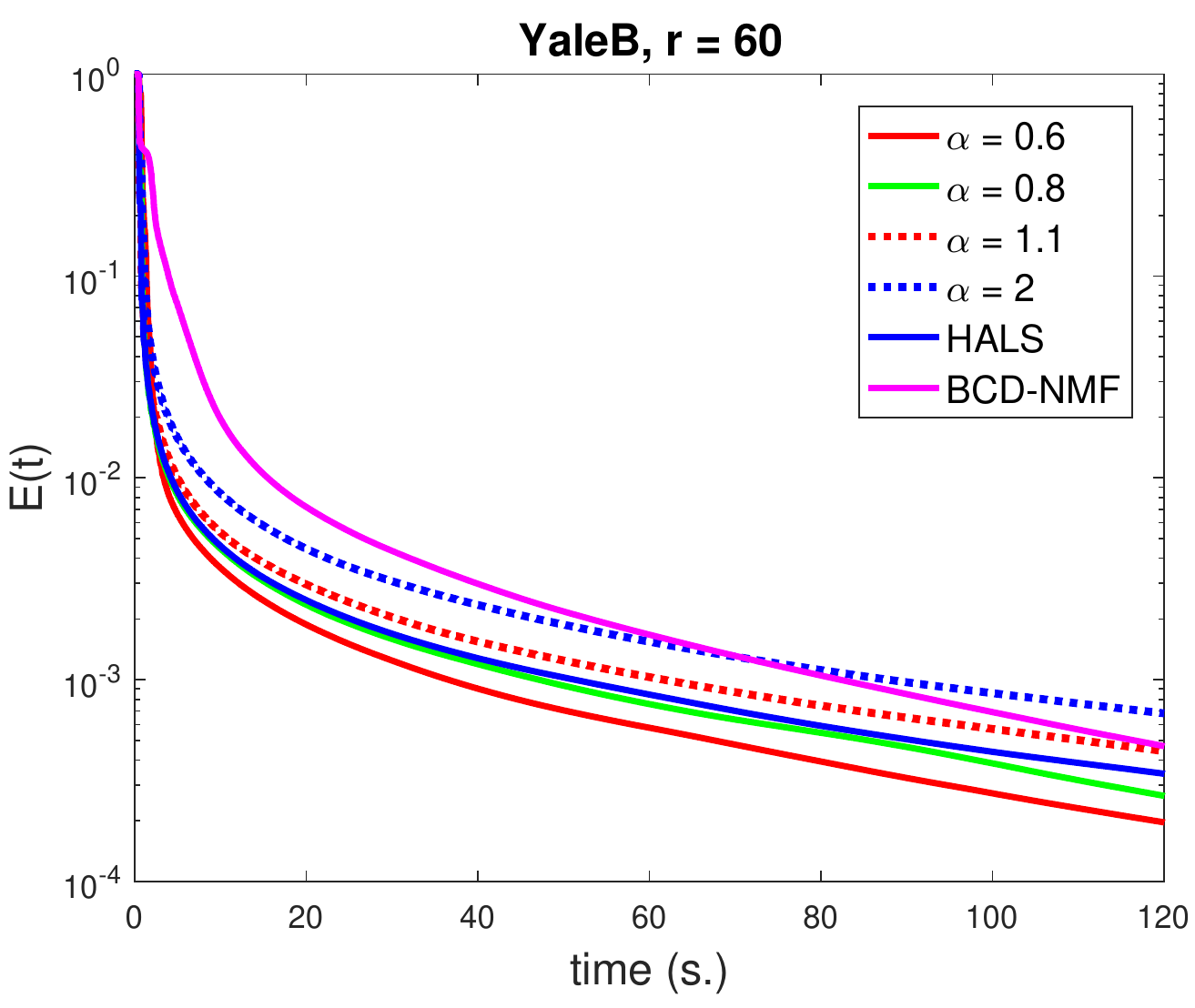}}
\caption{Average $E(t)$ of 30 independent trials for NMF on face datasets.}\label{Et_face}
\end{figure}

\begin{figure}[ht]
\centering
\subfigure[$\mathrm{T}^{\max}=3$]{\includegraphics[width=7cm]{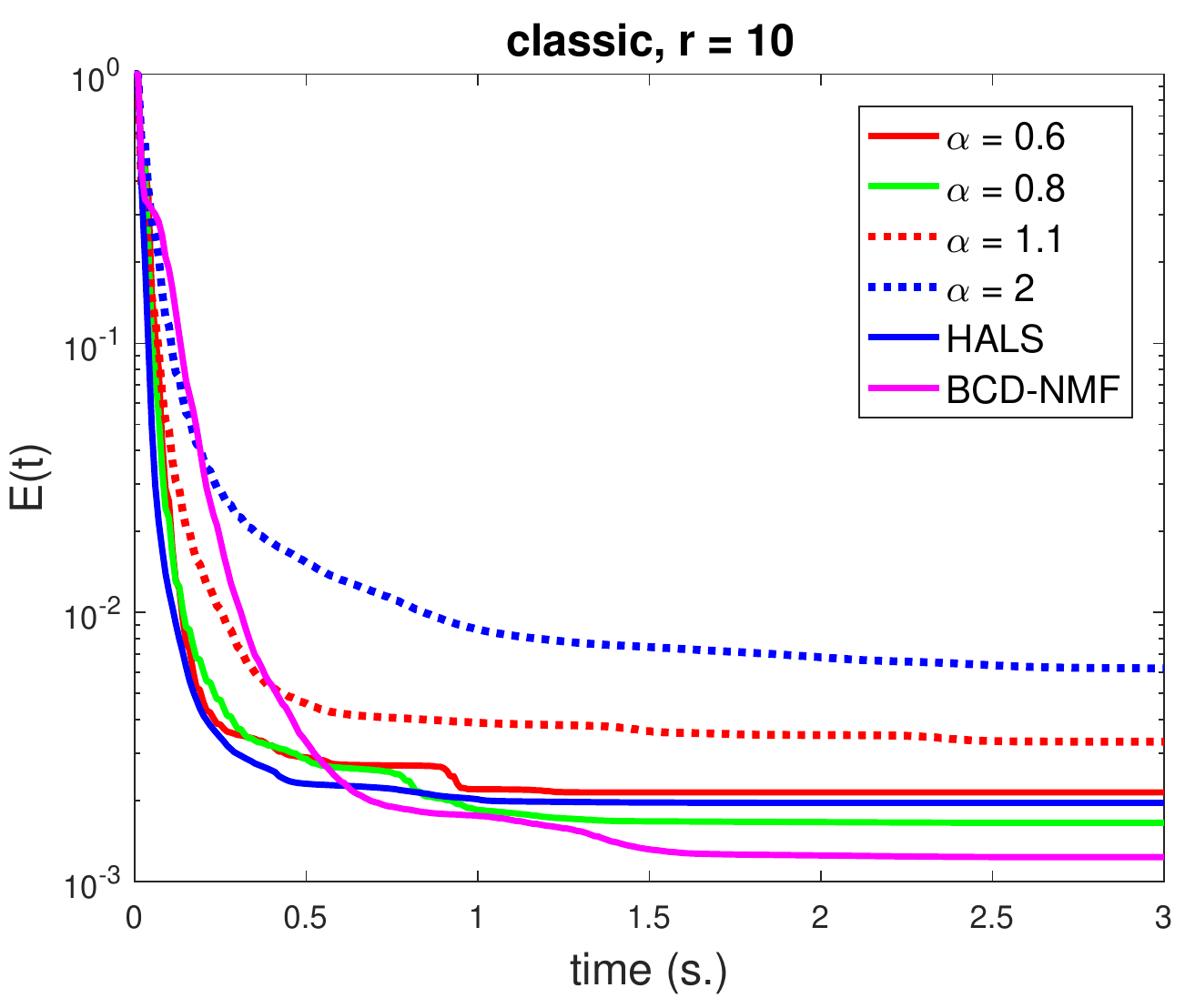}}
\subfigure[$\mathrm{T}^{\max}=6$]{\includegraphics[width=7cm]{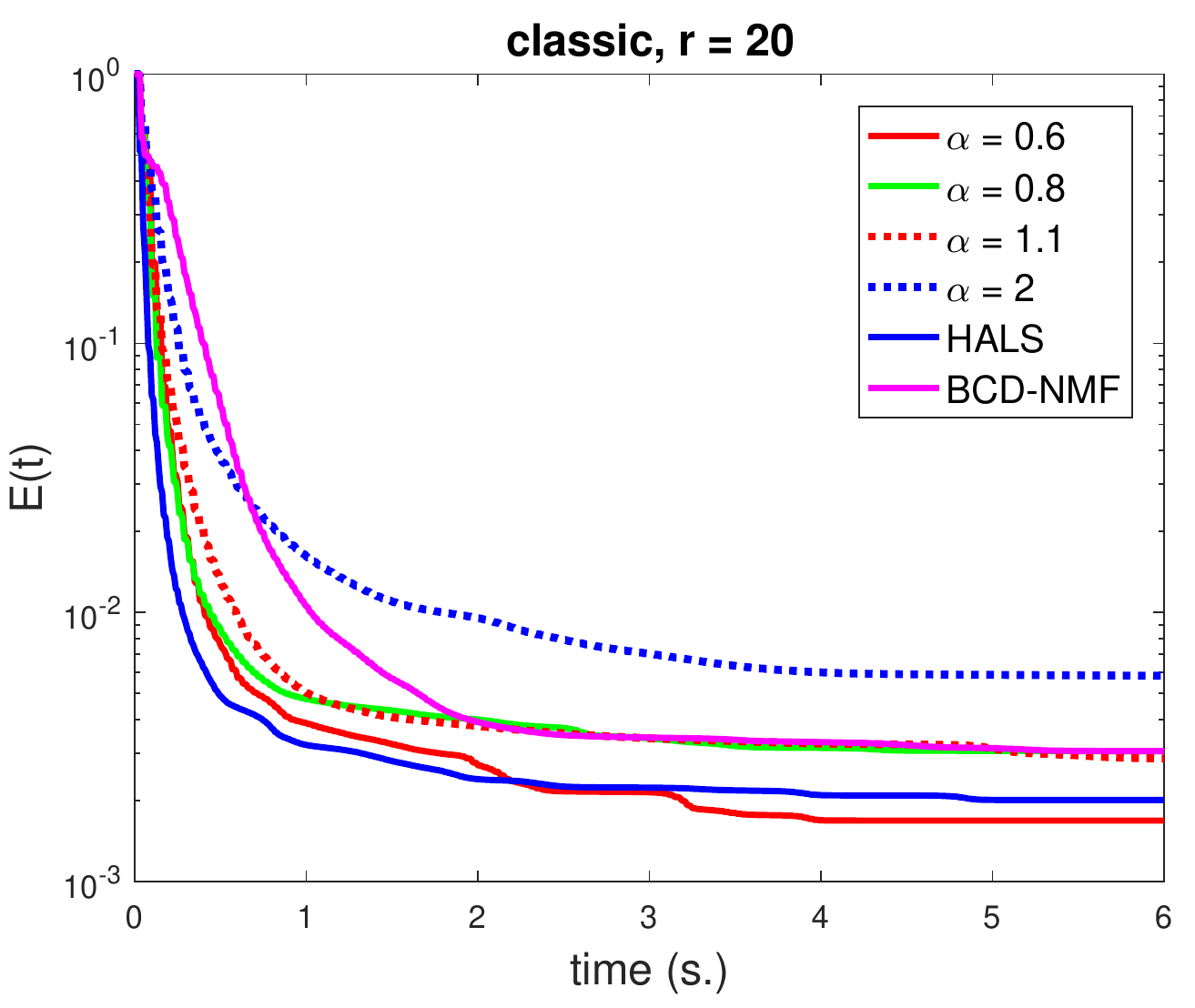}} \\
\subfigure[$\mathrm{T}^{\max}=3$]{\includegraphics[width=7cm]{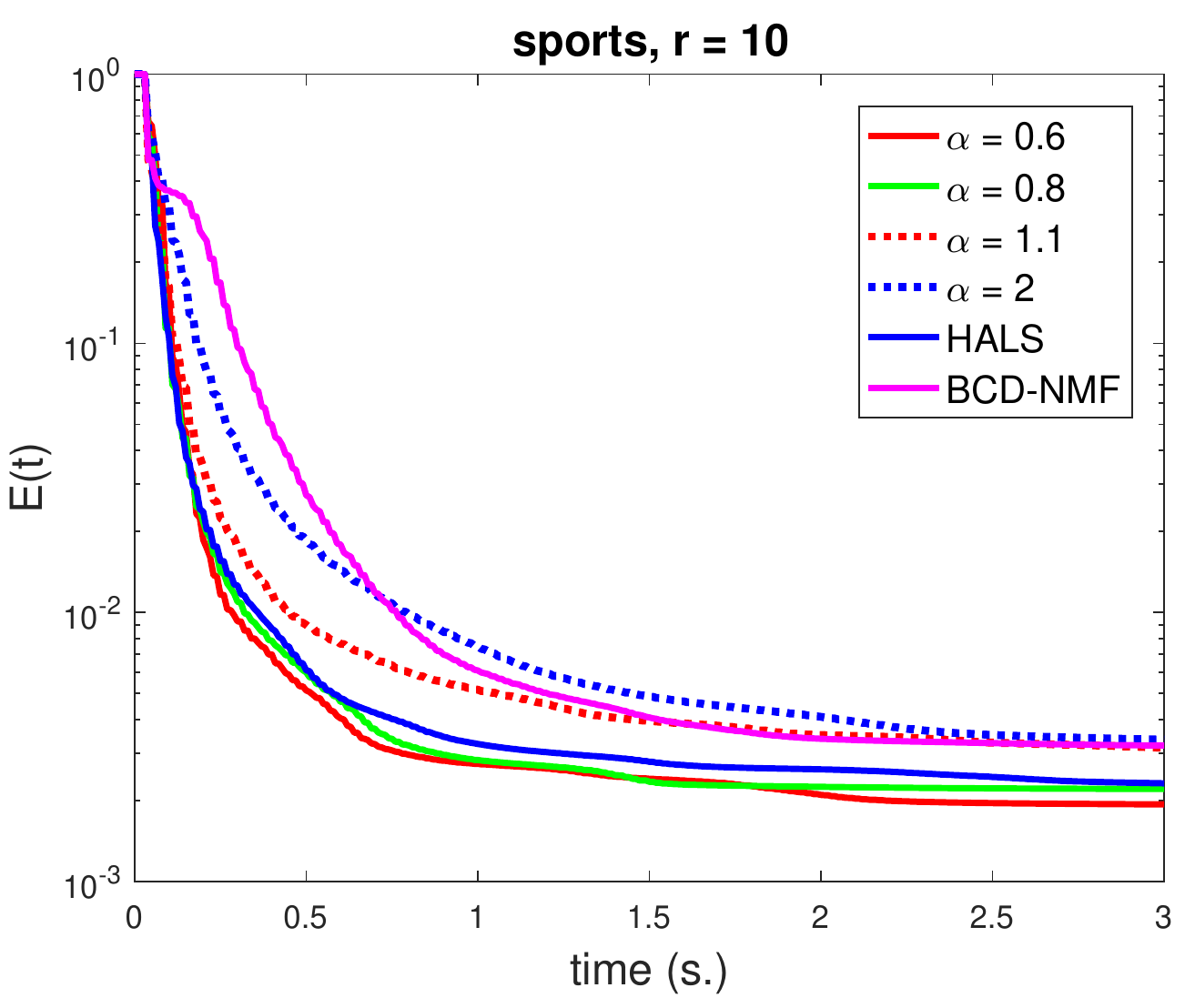}}
\subfigure[$\mathrm{T}^{\max}=6$]{\includegraphics[width=7cm]{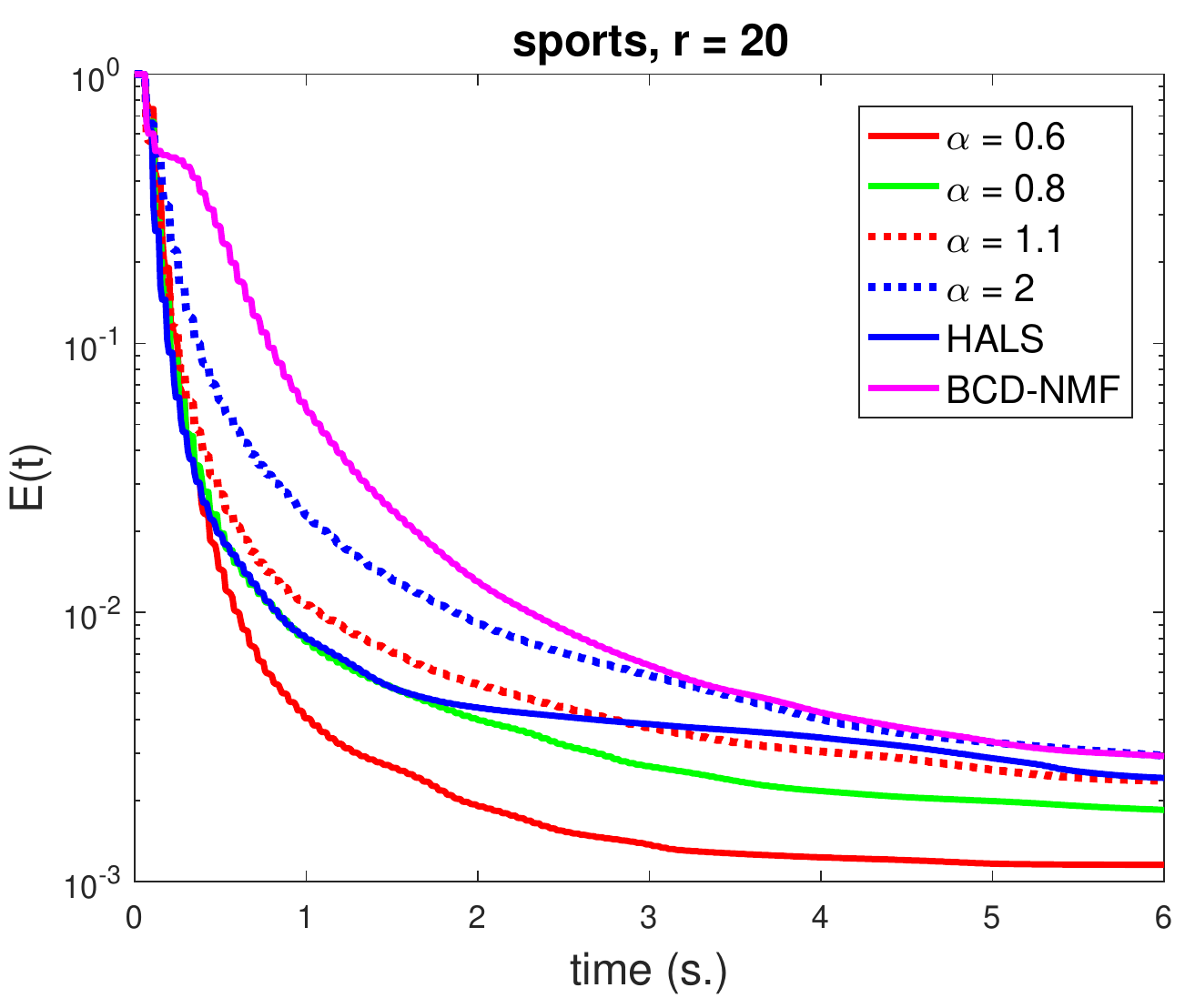}}  \\
\subfigure[$\mathrm{T}^{\max}=3$]{\includegraphics[width=7cm]{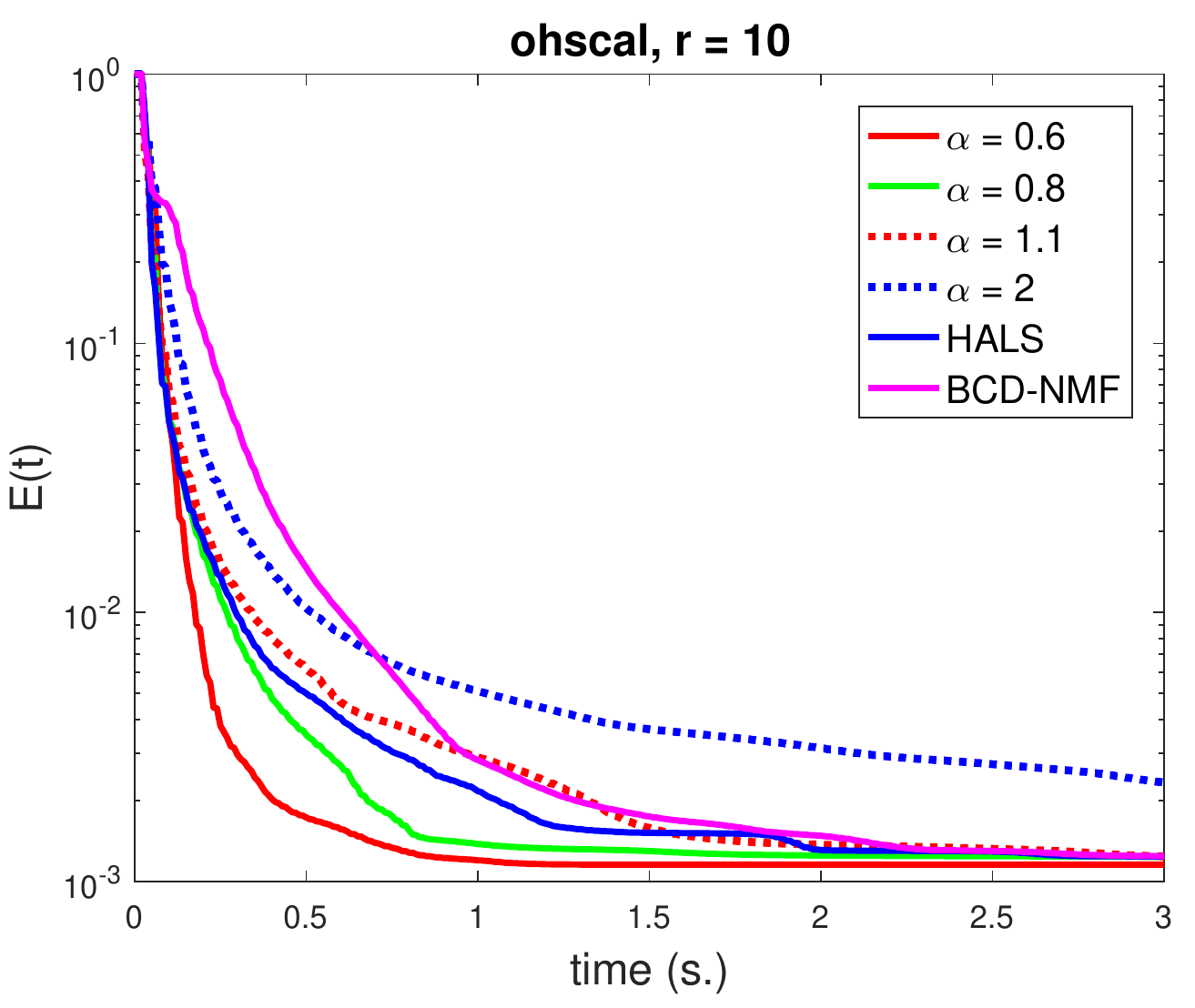}}
\subfigure[$\mathrm{T}^{\max}=6$]{\includegraphics[width=7cm]{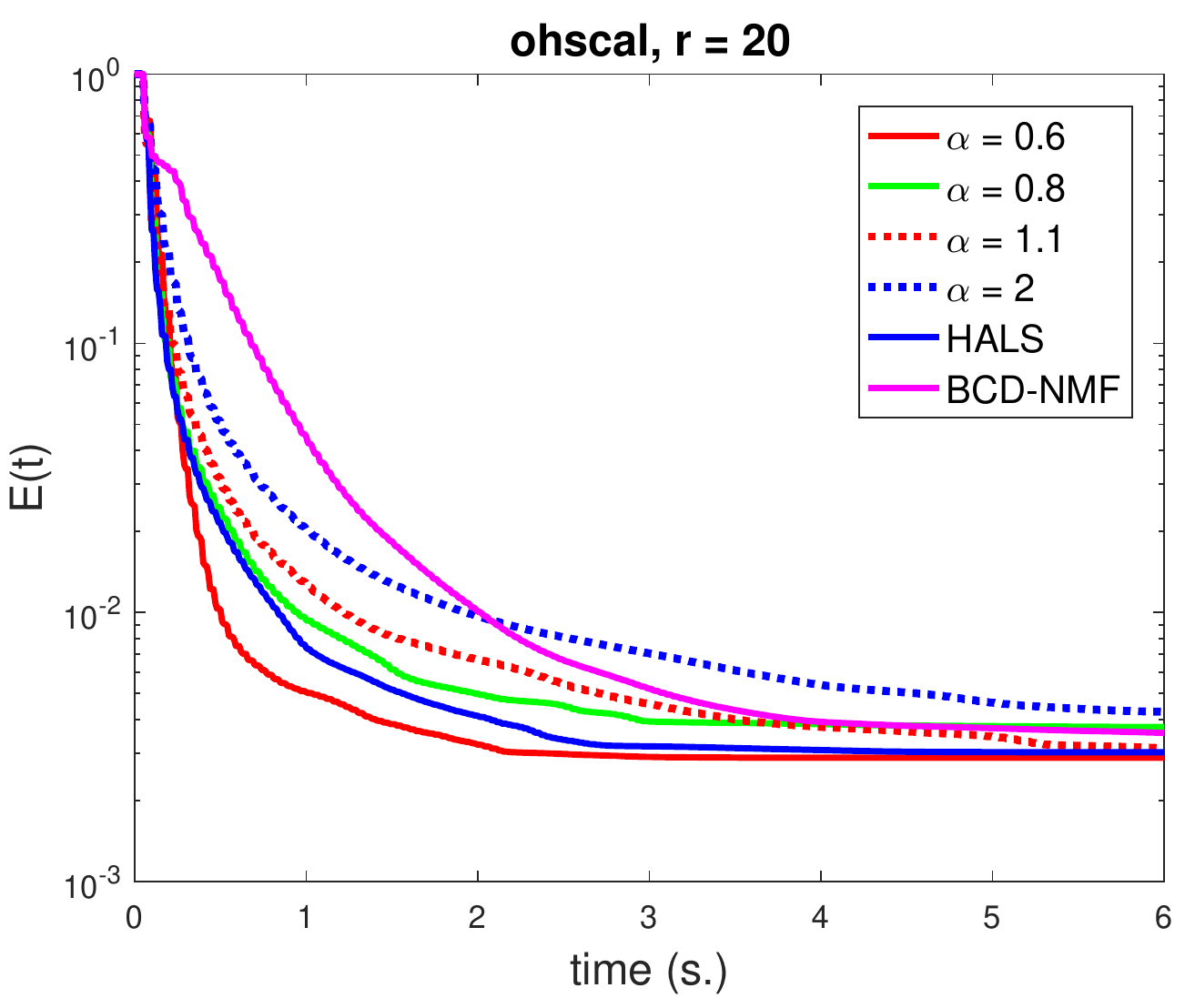}}
\caption{Average $E(t)$ of 30 independent trials for NMF on text datasets.}\label{Et_text}
\end{figure}

In our experiments, we initialize all the algorithms at the same random initial point $(X^0, \,Y^0)$ and terminate them \textit{only} by the maximum running time $\mathrm{T}^{\max}$. The specific values of $\mathrm{T}^{\max}$ are given in Fig.\,\ref{Et_face} and Fig.\,\ref{Et_text}. Additionally, we use the default settings for BCD-NMF. For NAUM, we choose $\alpha=0.6, 0.8, 1.1, 2$. We then plot the average $E(t)$ for each algorithm within time $\mathrm{T}^{\max}$.

Fig.\,\ref{Et_face} and Fig.\,\ref{Et_text} show the average $E(t)$ of 30 independent trials for NMF on face datasets and text datasets, respectively. From the results, we can see that NAUM with $\alpha=0.6$ performs best in most cases, and NAUM with $\alpha=0.6$ or $0.8$ always performs better than NAUM with $\alpha>1$. This shows that choosing $\alpha$ and $\beta$ under the weaker condition $\frac{1}{\alpha} +\frac{1}{\beta}=1$ (hence $\alpha$ can be small than 1) can improve the numerical performance of NAUM.

\subsection{Matrix completion}

We next consider a recent model for MC:
\begin{eqnarray}\label{mcmodel}
\min \limits_{X,Y}~~\frac{\eta}{2}\|X\|_* + \frac{\eta}{2}\|Y\|_*+\frac{1}{2}\left\|\mathcal{P}_{\Omega}(XY^{\top}-M)\right\|_F^2,
\end{eqnarray}
where $\eta>0$ is a penalty parameter, $\Omega$ is the index set of the known entries of $M$, and $\mathcal{P}_{\Omega}(Z)$ keeps the entries of $Z$ in $\Omega$ and sets the remaining ones to zero. This model was first considered in \cite{slc2016scalable,slc2016tractable} and was shown to be equivalent to Schatten-$\frac{1}{2}$ quasi-norm minimization. Encouraging numerical performance of this model has also been reported in \cite{slc2016scalable,slc2016tractable}. Note that \eqref{mcmodel} corresponds to \eqref{MFPmodel} with $\Psi(X)=\frac{\eta}{2}\|X\|_*$, $\Phi(Y)=\frac{\eta}{2}\|Y\|_*$ and $\mathcal{A}=\mathcal{P}_{\Omega}$. Thus, we can apply NAUM with \eqref{Xlinupdate} and \eqref{Ylinupdate} to solving \eqref{mcmodel}. The updates of $Z^k$, $U$ and $V$ are
\begin{eqnarray*}
\begin{aligned}
Z^k &= X^k(Y^k)^{\top} + {\textstyle\frac{\beta}{\alpha+\beta}}\mathcal{P}_{\Omega}\left(M-X^k(Y^k)^{\top}\right),     \\
U &= \mathcal{S}_{\eta/(2\mu_k)}\left(X^k-{\textstyle\frac{\alpha}{\mu_k}}(X^k(Y^k)^{\top}-Z^k)Y^k\right),             \\
V &= \mathcal{S}_{\eta/(2\sigma_k)}\left(Y^k-{\textstyle\frac{\alpha}{\sigma_k}}(U(Y^k)^{\top}-Z^k)^{\top}U\right).    \\
\end{aligned}
\end{eqnarray*}
Substituting $Z^k$ into $U$ and $V$ and using $\frac{1}{\alpha}+\frac{1}{\beta}=1$ gives
\begin{eqnarray}\label{substz2}
\begin{aligned}
&U = \mathcal{S}_{\eta/(2\mu_k)}\left(X^k\!-\!{\textstyle\frac{1}{\mu_k}}\left[\mathcal{P}_{\Omega}(X^k(Y^k)^{\top}\!-\!M)\right]Y^k\right),             \\
&V = \mathcal{S}_{\eta/(2\sigma_k)}\left(Y^k\!-\!{\textstyle\frac{\alpha}{\sigma_k}}Y^k(U\!-\!X^k)^{\top}U
\!-\!{\textstyle\frac{1}{\sigma_k}}\left[\mathcal{P}_{\Omega}(X^k(Y^k)^{\top}\!-\!M)\right]^{\top}U\right).    \\
\end{aligned}
\end{eqnarray}
Thus, similar to NAUM for NMF, we do not need to update $Z^k$ explicitly for MC.

We compare NAUM with proximal alternating linearized minimization (PALM), which was proposed in \cite{bst2014proximal} and was used to solve \eqref{mcmodel} in \cite{slc2016scalable,slc2016tractable}. For ease of future reference, we recall that the PALM for solving \eqref{mcmodel} is given by
\begin{eqnarray*}
\begin{aligned}
&X^{k+1} = \mathcal{S}_{\frac{\eta}{2\|Y^k\|^2}}\left(X^k-{\textstyle\frac{1}{\|Y^k\|^2}}\,
\left[\mathcal{P}_{\Omega}(X^k(Y^k)^{\top}-M)\right]Y^k\right), \\
&Y^{k+1} = \mathcal{S}_{\frac{\eta}{2\|X^{k+1}\|^2}}\left(Y^k-{\textstyle\frac{1}{\|X^{k+1}\|^2}}\,
\left[\mathcal{P}_{\Omega}(X^{k+1}(Y^k)^{\top}-M)\right]^{\top}X^{k+1}\right).
\end{aligned}
\end{eqnarray*}
For NAUM, we use the same parameter settings as in Section \ref{secNMF}, but choose $\alpha=0.4,0.6,1.1$. All the algorithms are initialized at the same random initialization $(X^0, \,Y^0)$\footnote{We use the Matlab commands: \texttt{X0 = randn(m, r); Y0 = randn(n, r);}} and terminated if one of the following stopping criteria is satisfied:
\begin{itemize}[leftmargin=6mm]
\item $\frac{|\mathcal{F}_{\mathrm{mc}}^{k} - \mathcal{F}_{\mathrm{mc}}^{k-1}|}{\mathcal{F}_{\mathrm{mc}}^{k}+1} \leq 10^{-4}$ holds for 3 consecutive iterations;

\item $\frac{\|X^k-X^{k-1}\|_F+\|Y^k-Y^{k-1}\|_F}{\|X^k\|_F+\|Y^k\|_F+1} \leq 10^{-4}$ holds;

\item the running time is more than 300 seconds,
\end{itemize}
where $\mathcal{F}_{\mathrm{mc}}^k:=\frac{\eta}{2}\|X^k\|_* + \frac{\eta}{2}\|Y^k\|_*+\frac{1}{2}\left\|\mathcal{P}_{\Omega}(X^k(Y^k)^{\top}-M)\right\|_F^2$ denotes the objective function value obtained by each algorithm at $(X^k, Y^k)$.

Table \ref{ResMc1} presents the numerical results of different algorithms for different problems, where two face datasets (CBCL and ORL) are used as our test matrices $M$ and a subset $\Omega$ of entries is sampled uniformly at random. In the table, $sr$ denotes the sampling ratio, i.e., a subset $\Omega$ of (rounded) $mn*sr$ entries is sampled; $r$ denotes the rank used for test; ``iter" denotes the number of iterations; ``Normalized fval" denotes the normalized function value $\frac{\mathcal{F}(X^*, \,Y^*) - \mathcal{F}_{\min}}{\mathcal{F}_{\max} - \mathcal{F}_{\min}}$, where $(X^*, \,Y^*)$ is obtained by each algorithm, $\mathcal{F}(X^*, \,Y^*)$ is the function value at $(X^*, \,Y^*)$ for each algorithm and $\mathcal{F}_{\max}$ (resp. $\mathcal{F}_{\min}$) denotes the maximum (resp. minimum) of the terminating function values obtained from \textit{all} algorithms in \textit{a} trial (one random initialization and $\Omega$); ``RecErr" denotes the recovery error $\frac{\|X^*(Y^*)^{\top}-M\|_F}{\|M\|_F}$. All the results presented are the average of 10 independent trials.

From Table \ref{ResMc1}, we can see that NAUM with $\alpha=0.4$ gives the smallest function values and the smallest recovery error within least CPU time in most cases. Moreover, NAUM with $\alpha=0.6$ also performs better than NAUM with $\alpha=1.1$ and PALM with respect to the function value and the recovery error in most cases. This again shows that a flexible choice of $\alpha$ and $\beta$ can lead to better numerical performances and the choice of $\alpha=0.4$ performs best for MC from our experiments.

\begin{table}[ht]
\setlength{\belowcaptionskip}{6pt}
\caption{Numerical results for MC on face datasets}\label{ResMc1}
\centering \tabcolsep 3pt
{\small
\begin{tabular}{|c|ccc|cccc|cccc|}
\hline
\multicolumn{1}{|c}{$\eta$} & \footnotesize{data} & $sr$ & $r$ & $\alpha=0.4$ & $\alpha=0.6$ & $\alpha=1.1$ & \footnotesize{PALM} & $\alpha=0.4$ & $\alpha=0.6$ & $\alpha=1.1$ & \footnotesize{PALM}\\
\hline \hline
\multicolumn{4}{|c|}{} & \multicolumn{4}{c|}{\footnotesize{iter}} & \multicolumn{4}{c|}{\footnotesize{Normalized fval}}\\
\hline
\multirow{8}{*}{5} & \multirow{4}{*}{\footnotesize{CBCL}}
& 0.5&30&  780 & 1189 & 3320 & 3306   &   1.13e-01 & {\bf7.50e-02} & 4.52e-01 & 1 \\
&&0.5&60&  921 & 1218 & 3850 & 4654   &   {\bf3.24e-02} & 5.10e-02 & 3.85e-01 & 1 \\
&&0.2&30&  1174 & 2366 & 4767 & 3573   &  {\bf8.01e-03} & 2.21e-01 & 6.87e-01 & 9.60e-01 \\
&&0.2&60&  1577 & 1919 & 5360 & 5037   &  {\bf1.03e-02} & 8.95e-02 & 8.08e-01 & 8.86e-01 \\
\cline{2-4}
&\multirow{4}{*}{\footnotesize{ORL}}
& 0.5&30&  1218 & 1243 & 1241 & 1468   &  {\bf0} & 2.94e-01 & 5.06e-01 & 1 \\
&&0.5&60&  1049 & 1051 & 1051 & 1327   &  {\bf0} & 1 & 4.00e-01 & 7.73e-01 \\
&&0.2&30&  2074 & 325 & 385 & 2691     &  {\bf2.59e-03} & 7.01e-01 & 1 & 1.31e-01 \\
&&0.2&60&  1551 & 1551 & 356 & 2222    &  {\bf0} & 3.82e-01 & 1 & 2.12e-01 \\
\hline
\multirow{8}{*}{10} & \multirow{4}{*}{\footnotesize{CBCL}}
& 0.5&30&  457 & 654 & 1793 & 1935   &   {\bf2.20e-02} & 1.29e-01 & 3.60e-01 & 9.81e-01 \\
&&0.5&60&  514 & 594 & 1950 & 2559   &   2.65e-01 & {\bf1.15e-01} & 3.79e-01 & 8.71e-01 \\
&&0.2&30&  627 & 1313 & 2513 & 2116   &  {\bf1.91e-02} & 3.75e-02 & 8.35e-01 & 7.79e-01 \\
&&0.2&60&  866 & 1095 & 2713 & 2889   &  {\bf2.07e-02} & 2.89e-02 & 9.22e-01 & 4.86e-01 \\
\cline{2-4}
&\multirow{4}{*}{\footnotesize{ORL}}
& 0.5&30&  1003 & 1186 & 1192 & 1402   & {\bf3.30e-02} & 1.47e-01 & 4.30e-01 & 1 \\
&&0.5&60&  975 & 1009 & 1012 & 1276   &  {\bf0} & 8.58e-01 & 6.11e-01 & 9.99e-01 \\
&&0.2&30&  1409 & 364 & 411 & 2646   &   {\bf0} & 7.16e-01 & 1 & 8.10e-02 \\
&&0.2&60&  1241 & 1504 & 376 & 2185   &  {\bf4.05e-06} & 3.97e-02 & 1 & 2.21e-01 \\
\hline \hline

\multicolumn{4}{|c|}{} & \multicolumn{4}{c|}{\footnotesize{CPU time}} & \multicolumn{4}{c|}{\footnotesize{RecErr}} \\
\hline
\multirow{8}{*}{5} & \multirow{4}{*}{\footnotesize{CBCL}}
& 0.5&30& {\bf35.56} & 54.14 & 151.23 & 119.05   &   {\bf1.05e-01} & {\bf1.05e-01} & 1.06e-01 & 1.08e-01 \\
&&0.5&60& {\bf57.66} & 76.09 & 240.19 & 206.47   &   {\bf8.81e-02} & 9.02e-02 & 9.04e-02 & 8.99e-02 \\
&&0.2&30& {\bf34.04} & 68.57 & 137.97 & 75.56    &   {\bf1.37e-01} & {\bf1.37e-01} & 1.38e-01 & 1.43e-01 \\
&&0.2&60& {\bf72.01} & 87.82 & 245.21 & 147.08   &   {\bf1.34e-01} & 1.35e-01 & 1.35e-01 & 1.36e-01 \\
\cline{2-4}
&\multirow{4}{*}{\footnotesize{ORL}}
& 0.5&30& {\bf294.20} & 300 & 300 & 300   &   {\bf1.72e-01} & 1.84e-01 & 2.01e-01 & 2.12e-01 \\
&&0.5&60& 300 & 300 & 300 & 300           &   {\bf1.66e-01} & 2.11e-01 & 2.05e-01 & 2.11e-01 \\
&&0.2&30& 300 & {\bf47.35} & 55.86 & 300  &   {\bf2.08e-01} & 3.04e-01 & 3.81e-01 & 2.24e-01 \\
&&0.2&60& 300 & 300 & {\bf69.21}  & 300   &   {\bf2.16e-01} & 2.35e-01 & 3.49e-01 & 2.61e-01 \\
\hline
\multirow{8}{*}{10} & \multirow{4}{*}{\footnotesize{CBCL}}
& 0.5&30& {\bf21.01} & 30.12 & 82.45 & 70.32   &   {\bf1.16e-01} & 1.19e-01 & 1.18e-01 & 1.17e-01 \\
&&0.5&60& {\bf32.40} & 37.38 & 122.51 & 113.80 &   {\bf1.09e-01} & 1.11e-01 & 1.14e-01 & 1.11e-01 \\
&&0.2&30& {\bf18.15} & 38.01 & 72.84 & 44.62   &   {\bf1.60e-01} & 1.61e-01 & 1.62e-01 & {\bf1.60e-01} \\
&&0.2&60& {\bf39.13} & 49.37 & 123.74 & 83.52  &   1.57e-01 & 1.57e-01 & 1.58e-01 & {\bf1.56e-01} \\
\cline{2-4}
&\multirow{4}{*}{\footnotesize{ORL}}
& 0.5&30& {\bf252.15} & 300 & 300 & 300        &   {\bf1.71e-01} & 1.77e-01 & 1.95e-01 & 2.08e-01 \\
&&0.5&60& {\bf289.57} & 300 & 300 & 300        &   {\bf1.53e-01} & 2.01e-01 & 2.03e-01 & 2.09e-01 \\
&&0.2&30& 207.22 & {\bf53.08} & 60.54 & 300    &   {\bf1.95e-01} & 3.06e-01 & 3.83e-01 & 2.14e-01 \\
&&0.2&60& 243.45 & 295.60 & {\bf74.09} & 300   &   {\bf1.87e-01} & 1.95e-01 & 3.60e-01 & 2.36e-01 \\
\hline
\end{tabular}}
\end{table}

\section{Concluding remarks}\label{secconc}
In this paper, we consider a class of matrix factorization problems involving two blocks of variables. To solve this kind of possibly nonconvex, nonsmooth and non-Lipschitz problems, we introduce a specially constructed potential function $\Theta_{\alpha,\beta}$ defined in \eqref{defpofun} which contains one auxiliary block of variables. We then develop a non-monotone alternating updating method with a suitable line search criterion based on this potential function. Unlike other existing methods such as those based on alternating minimization, our method essentially updates the two blocks of variables alternately by solving subproblems related to $\Theta_{\alpha,\beta}$ and then updates the auxiliary block of variables by an explicit formula (see \eqref{Zupdate}). Using the special structure of $\Theta_{\alpha,\beta}$, we demonstrate how some efficient computational strategies for NMF can be used to solve the associated subproblems in our method. Moreover, under some mild conditions, we establish that the sequence generated by our method is bounded and any cluster point of the sequence gives a stationary point of our problem. Finally, we conduct some numerical experiments for NMF and MC on real datasets to illustrate the efficiency of our method.

Note that the parameter $\alpha$ (and $\beta=\alpha/(\alpha-1)$) plays a significant role in our NAUM. Although it has been observed in our experiments that a relatively small $\alpha$ (e.g., 0.6, 0.8) can improve the numerical performance of NAUM, how to choose an optimal $\alpha$ is still unknown. In view of the recent work \cite{op2017adaptive} on adaptively choosing the extrapolation parameter in FISTA for solving a class of possibly nonconvex problems, it may be possible to derive a strategy to adaptively update $\alpha$ in our NAUM. This is a possible future research topic.

\section*{Acknowledgments}

\noindent The authors are grateful to Nicolas Gillis for his helpful comments. The authors are also grateful to the editor and the anonymous referees for their valuable suggestions and comments, which helped improve this paper.


\end{document}